\documentclass[journal]{IEEEtran}
\pdfoutput=1

\usepackage{graphicx}
\graphicspath{ {Pictures/} }
\ifCLASSINFOpdf
\else
\fi
\usepackage{amsmath,amsfonts,amsthm} 
\allowdisplaybreaks
\usepackage{enumerate}
\usepackage{eurosym}
\usepackage{booktabs}
\usepackage{accents}
\usepackage{color}
\allowdisplaybreaks

\begin{document}

\title{Enabling Active/Passive Electricity Trading in Dual-Price Balancing Markets}
	
\author{Nicol\`o~Mazzi,\IEEEmembership{}
		Alessio~Trivella,\IEEEmembership{}
		and~Juan~M.~Morales,~\IEEEmembership{Senior Member,~IEEE}
		\thanks{N. Mazzi is with the School of Mathematics, University of Edinburgh, Edinburgh,
			UK (e-mail: nicolo.mazzi@ed.ac.uk), and is supported by the Engineering and Physical Sciences Research Council (EPSRC) through the CESI project (EP/P001173/1).}
		\thanks{A. Trivella is with DTU Management Engineering, Technical University of Denmark, Kgs. Lyngby,
			Denmark (e-mail: atri@dtu.dk), and is supported by the Innovation Fund Denmark through the project SAVE-E (4106-00009B).}
		\thanks{J. M. Morales is with the Department of Applied Mathematics, University of M\'{a}laga, M\'{a}laga, Spain (e-mail: juan.morales@uma.es). The work of Juan M. Morales is partly funded by the European Research Council (ERC) under the European Union’s Horizon 2020 research and innovation programme (grant agreement No. 755705), by the Spanish Research Agency through project ENE2017-83775-P and by the Research Funding Program for Young Talented Researchers of the University of M\'{a}laga through project PPIT-UMA-B1-2017/18.}
	}
	
	\maketitle

	\begin{abstract}
		In electricity markets with a dual-pricing scheme for balancing energy, controllable production units typically participate in the balancing market as ``active'' actors by offering regulating energy to the system, while renewable stochastic units are treated as ``passive'' participants that create imbalances and are subject to less competitive prices. Against this background, we propose an innovative market framework whereby the participant in the balancing market is allowed to act as an active agent (i.e., a provider of regulating energy) in some trading intervals and as a passive agent (i.e., a user of regulating energy) in some others. To illustrate and evaluate the proposed market framework, we consider the case of a virtual power plant (VPP) that trades in a two-settlement electricity market composed of a day-ahead and a dual-price balancing market. We formulate the optimal market offering problem of the VPP as a three-stage stochastic program, where uncertainty is in the day-ahead electricity prices, balancing prices and the power output from the renewable units.
		Computational experiments show that the VPP expected revenues can increase substantially compared to an active-only or passive-only participation, and in the paper we discuss how the variability of the stochastic sources affects the balancing market participation choice.
	\end{abstract}

	
	\IEEEpeerreviewmaketitle
	
	\section*{Nomenclature}
	\addcontentsline{toc}{section}{Nomenclature}
	\subsection*{Indices and Sets} \vspace{0pt}
	\begin{IEEEdescription}[\IEEEusemathlabelsep\IEEEsetlabelwidth{$ \text{Beta}(\quad) $}]
		\item [$ i,i' \in I $]    Indices of day-ahead market price scenarios
		\item [$ j,j'  \in J $]    Indices of balancing market price scenarios
		\item [$ \omega \in W $]   Index of renewable energy generation scenarios
		\item [$ k \in K $]   Index of time intervals
		\item [$ \Pi^{\mathrm{DA}} $]    Feasible region of the day-ahead market offers
		\item [$ \Pi^{\mathrm{Act}} $]   Feasible region of the active participation
		\item [$ \Pi^{\mathrm{Pas}} $]   Feasible region of the passive participation
		\item [$ \Omega $]    Feasible region of the VPP's operation
		\item [$ \Gamma $]    Feasible region of \textit{Active/Passive} participation
	\end{IEEEdescription}
	
	\subsection*{Parameters} \vspace{0pt}
	\begin{IEEEdescription}[\IEEEusemathlabelsep\IEEEsetlabelwidth{$ \text{Beta}(\quad) $}]
		\item [$ \lambda^{\textrm{DA}}_{ik} $] Day-ahead market price (\euro/MWh)
		\item [$ \lambda^{\textrm{BA}}_{ijk} $] Balancing market price (\euro/MWh)
		\item [$ E_{\omega k} $] Wind (or solar) power generation (MWh)
		\item [$ \overline{E} $] Capacity of the wind (or solar) unit (MW)
		\item [$ \overline{D} $] Capacity of the thermal unit (MW)
		\item [$ \underline{D} $] Minimum power limit of the thermal unit (MW)
		\item [$ R^{\textrm{UP}},R^{\textrm{DW}} $] Thermal unit ramp-up and -down limits (MW/h)
		\item [$ \overline{P}^{(\uparrow)},\overline{P}^{(\downarrow)} $] Charging/discharging power limits (MW)
		\item [$ \underline{L},\overline{L}$] Minimum/maximum level of the storage (MWh)
		\item [$ \eta $] Round-trip efficiency of the energy storage 
		\item [$ C $] Marginal cost of the thermal unit (\euro/MWh)
		\item [$ C_0 $] Fixed cost of the thermal unit (\euro)
		\item [$ \pi^{\textrm{DA}}_i$] Probability of day-ahead price scenario $i$
		\item [$ \pi^{\textrm{BA}}_{ij} $] Probability of balancing price scenario $j$, provided that day-ahead price scenario $i$ realizes
		\item [$ \pi^{\textrm{E}}_{\omega} $]  Probability of renewable energy production scenario $\omega$
	\end{IEEEdescription}
	
	\subsection*{Variables} \vspace{0pt}
	\begin{IEEEdescription}[\IEEEusemathlabelsep\IEEEsetlabelwidth{$ \text{Beta}(\quad) $}]
		\item [$ q^{\textrm{DA}}_{ik} $]  Quantity offer at day-ahead market (MWh)
		\item [$ q^{\textrm{UP}}_{ijk}, q^{\textrm{DW}}_{ijk} $]    Up/down regulation quantity offer (MWh)
		\item [$ q^{(+)}_{i\omega k}, q^{(-)}_{i\omega k} $]    Positive/negative real-time deviation (MWh)
		\item [$ d_{ij\omega k} $] Thermal unit energy production (MWh)
		\item [$ p^{(\uparrow)}_{ij\omega k},p^{(\downarrow)}_{ij\omega k} $] Charging/discharging quantities (MWh)
		\item [$ \ell_{ij\omega k} $] Energy storage level (MWh)
		\item [$ \hat{\rho}^{\textrm{DA}}_{k} $]    Expected day-ahead market revenue (\euro)
		\item [$ \hat{\rho}^{\textrm{Act}}_{k} $]   Expected revenue of the active participation (\euro)
		\item [$ \hat{\rho}^{\textrm{Pas}}_{k} $]   Expected revenue of the passive participation (\euro)
		\item [$ \hat{c}_{k} $] Expected operational cost (\euro)
		\item [$ u_{ij\omega k} $]   Commitment (binary) status of the thermal unit
		\item [$ \epsilon_{i k} $]   Auxiliary binary variables to enforce complementarity of the \textit{Active/Passive} participation
	\end{IEEEdescription}

\section{Introduction}\label{sec:Introduction}
	
	\IEEEPARstart{S}{ociety} is moving towards using more renewable energy sources to decrease the dependency on fossil fuels. Governments, seeking to increase the share of renewable energy, typically support stochastic power sources such as wind and solar power by means of subsidies. However, with the steep growth and decreasing cost of renewable energy generation experienced in the recent years, stochastic producers are increasingly required to be financially responsible for the imbalances created in the real-time. Accordingly, renewable energy producers access the balancing electricity market as ``passive'' actors by settling the deviation from the day-ahead contracted schedule at a less favorable power price.
	
	The imbalance of the system, often caused by forecast errors of power demand and renewable generation, is restored by rescheduling the market position of the ``active'' participants in the balancing market. Such producers offer to the Transmission System Operator (TSO) the flexibility to upward or downward adjust their day-ahead contracted schedule, provided to be remunerated at a more convenient price. However, to qualify as regulators in the balancing market, generators must fulfill specific requirements from the TSO, which include the ability to always meet the contracted schedule except for unpredictable unit failures. As a consequence, only conventional generators can currently be qualified as active balancing market participants. In this context, it is quite straightforward to distinguish between passive participants (i.e., stochastic producers) that regularly deviate from their contracted schedule, and active participants (i.e., conventional producers) that can consistently respect their market position and offer additional regulating energy to the TSO. 
	
	In this context, virtual power plants (VPPs), i.e., clusters of combined generating units, storage systems and flexible loads that act as a single participant in the electricity market \cite{morales2013integrating}, lie somehow in the middle between the two classes of stochastic and conventional generators. There is an increasing interest towards VPPs that is associated with the possibility of internally handling the forecast errors caused by renewable energy units. However, we believe that VPPs have more potential than solely self-balancing their stochastic energy production. Specifically, if VPPs were able to provide regulating energy to the TSO, when available, then they would bring a major benefit to the operation of the system. Indeed, a VPP often uses only a portion of its flexibility to balance the stochastic power production within the cluster, and additional flexibility could be used to compensate for the imbalance created by renewable energy generation outside the cluster. This potential is not exploited by the TSO due to the current balancing market setup, which enforces a clear distinction between active and passive participation in the balancing market. This results in classifying a VPP that includes stochastic generation as a passive participant, since it can hardly fulfill the qualification procedure for being an active regulator.
	
	In this paper, we propose an innovative and more flexible market framework where we loosen the TSO's requirements for producers to qualify as active participants in the balancing market. Specifically, our idea is to allow VPPs, and market participants in general, to actively sell regulating energy in some trading intervals and passively deviate from their schedule in others. This relaxes the duality between only-active and only-passive participation, and allows to exploit the full potential of VPPs in electricity markets. What we propose can be seen as a new ``market product'' or ``market bid format'', called \textit{Active/Passive} offers, tailored to the characteristics of VPPs. Note that market products that accommodate the needs of a specific technology are not uncommon. For instance, the ``linked block orders'' in the Nord Pool market are designed for conventional generators that want to include start-up and/or shut-down costs within their market offers. It is worth specifying that this new market product considers the potential reliability issues arising from breaking the distinction between active and passive actors. Indeed, VPPs that commit to sell regulating energy during a trading interval are prevented from creating an imbalance in the same interval.
	
	To illustrate and evaluate the potential impact of the proposed balancing market framework, we take the perspective of a VPP offering in a two-settlement electricity market that allows for \textit{Active/Passive} offers. Our goal is determining whether VPPs would gain from actively offering regulating energy, or would instead choose to behave as passive actors. 
    
    In the remainder of this section, we presents the literature review on optimal offering and modeling of VPPs in electricity markets in Section \ref{subsec:Literature_Review}. Section \ref{subsec:Approach_and_Contributions} presents our approach and the contributions of this work. Finally, Section  \ref{subsec:Paper_Organization} describes the structure of the paper.

\subsection{Literature Review}\label{subsec:Literature_Review}
	
	The problem of determining the optimal market offer for a stochastic power producer has been widely studied in the literature. In \cite{bremnes2004}, \cite{pinson2007trading}, and \cite{bitar2012bringing}, the optimal quantity to be submitted in the day-ahead market is derived as a quantile of the probability distribution of the future wind or solar power production. This idea is extended in \cite{dent2011opportunity} including the correlation between future wind production and real-time prices. In \cite{morales2010short}, the optimal offering strategy for a wind power producer is solved using stochastic programming.
	These models consider the stochastic producer as a passive actor in the balancing market, i.e., the balancing stage is only used to settle deviations from the day-ahead contracted schedule.
	
	Similarly, several models have been developed to derive the optimal offering strategy for a conventional production unit. In \cite{Arroyo2000}, \cite{Arroyo2004}, and \cite{Conejo2004}, the feasible operating region of a thermal unit is formulated using a mixed-integer linear program (MILP). Other papers (e.g., \cite{Conejo2002}, \cite{Baillo2004}, \cite{conejo2010decision}, and \cite{Maenhoudt2014}) combine such operating region with an electricity market trading problem obtaining different offering strategies. In contrast to stochastic units, conventional units are modeled as active participants in the balancing market, i.e., they access the balancing market to offer regulating energy to the system operator.
	
	The optimal participation of a VPP in an electricity market has been less investigated. In \cite{ruiz2009direct}, a direct load control algorithm is used for managing an aggregate of controllable loads. A deterministic offering model for a cluster consisting of combined heat and power units and a wind farm participating in a dual-price balancing market is given in \cite{hellmers2016operational}. Similarly, \cite{zapata2014comparative} presents a comparative (deterministic) study to reduce the real-time imbalances of a VPP composed of combined heat and power units and PV solar units. Papers \cite{mashhour2011biddingI} and \cite{mashhour2011biddingII} study the bidding problem of a VPP composed of dispatchable units in a joint market for energy and reserve. The authors consider a deterministic setting and formulate a mixed-integer non-linear program solved using a genetic algorithm. In \cite{peik2013decision}, a VPP offering strategy is formulated as a unit commitment problem where point estimates are used to model the uncertainty in market prices and power generation. Reference \cite{pandvzic2013virtual} proposes a stochastic MILP to derive the optimal self-scheduling of a VPP, considering a weekly time horizon and including long-term bilateral contracts and technical constraints of the units. Subsequently, the authors of \cite{pandvzic2013offering} develop a two-stage stochastic offering model to maximize the expected profit of a VPP with uncertainty in electricity prices and power production. Other works, e.g., \cite{kardakos2016optimal}, include the electricity market clearing process within the optimal offering strategy, resulting in a hierarchical stochastic optimization model. Finally, we refer to \cite{morales2013integrating} for a general VPP modeling approach in which different combinations of generating units, flexible loads, and storage systems are examined.
	
\subsection{Contributions}\label{subsec:Approach_and_Contributions}
	The contribution of our paper is threefold:
	\begin{itemize}
		\item[$i$] \underline{Conceptual contribution}. We propose a  new market product, or market bid format, called 		\textit{Active/Passive} offers, tailored to the characteristics of VPPs. This innovative market product has the potential to fully exploit the flexibility of the dispatchable generators of a VPP, i.e., balancing the forecast errors of the stochastic power units both within and outside the VPP cluster. Our proposed market product could also be adapted or extended to more ``sophisticated" products, e.g., the TSO may require the VPP to operate in the same mode (either active or passive) for at least a given number of trading intervals.
		\item[$ii$] \underline{Modeling contribution}. We develop the optimal offering strategy for a VPP in an electricity market that allows for the novel \textit{Active/Passive} offers. Our new market concept forces to rethink the VPP strategy and has a serious impact on the mathematical formulation. Indeed, the VPP offering strategies from the extant literature model the VPP either as a passive balancing market actor (e.g., \cite{ruiz2009direct,hellmers2016operational,zapata2014comparative,peik2013decision,pandvzic2013virtual,pandvzic2013offering} and \cite{kardakos2016optimal}) which solely settles the real-time deviations that are not self-balanced, or as an active actor (e.g., \cite{mashhour2011biddingI} and \cite{mashhour2011biddingII}) that can only sell regulating energy and is not allowed to deviate from the contracted schedule. In contrast, we account for the \textit{Active/Passive} participation by developing a three-stage stochastic program with recourse, which is more involved than the offering strategies from the literature (usually deterministic or two-stage stochastic programs). The \textit{Active/Passive} offering strategies also require additional binary variables and new constraints to enforce the complementarity between the active and the passive participation. This has an impact on the computational time needed to derive the operating strategy and should be thus accounted for when considering introducing the \textit{Active/Passive} offers in an electricity market.
		\item[$iii$] \underline{Practical impact}. In addition to proposing a new concept (i.e., \textit{Active/Passive} offers) and a new model (i.e., bidding strategy with \textit{Active/Passive} offers), we also validate numerically the interest of VPPs towards such offers. Our case studies show that (i) the VPP is often willing to actively sell regulating energy, and (ii) the \textit{Active/Passive} offer always outperforms the active-only and passive-only strategy in terms of expected revenues. This increased profitability may also encourage the aggregation of different technologies into a VPP.
	\end{itemize}
		
\subsection{Paper Structure}\label{subsec:Paper_Organization}
		
	The rest of this paper is organized as follows. We start in Section \ref{sec:Modelling_Assumption} by presenting the electricity market framework, the VPP structure, and characterizing the uncertainty. In Section \ref{sec:Optimal_Offering_Strategy}, we formulate the \textit{Active/Passive} offering strategy as a three-stage stochastic program. Section \ref{sec:Example} gives an explanatory example and Section \ref{sec:Case_Study} compares the \textit{Active/Passive} offering strategy with the \textit{Passive} and \textit{Active} strategies for different VPP configurations. Conclusions are drawn in Section \ref{sec:Conclusions}.

\section{Market Framework and Modeling Assumptions}\label{sec:Modelling_Assumption}
		
\subsection{Electricity Market Framework}\label{subsec:Electricity_Market_Framework}
		
	We consider a two-settlement electricity market composed of a day-ahead and a balancing market. The day-ahead market is cleared at noon for all 24 hourly trading intervals of the following day. The accepted day-ahead market offers are settled under a uniform pricing scheme. Subsequently, closer to the real-time operation, a separate balancing market is cleared for each hourly interval, one hour before operation. At the balancing stage the active participants submit their offers, in the form of non-decreasing offer curves, for the provision of regulating energy to the TSO. The accepted offers are priced under a uniform pricing scheme. 
    Specifically, if an active participant submits an up-regulation offer of quantity $q^{\mathrm{UP}}_k$ during interval $k$, this is accepted if the system actually requires up-regulation, and the offered price is lower than the balancing market price $\lambda^{\mathrm{BA}}_k$ (merit order). If these conditions are satisfied, then the active participant receives an income of $\lambda^{\mathrm{BA}}_k q^{\mathrm{UP}}_k$. Similarly, a down-regulation offer of quantity $q^{\mathrm{DW}}_k$ at $k$ is accepted if the system needs down-regulation, and the offered price is higher than $\lambda^{\mathrm{BA}}_k$, resulting in a cash flow of $-\lambda^{\mathrm{BA}}_k q^{\mathrm{DW}}_k$.

Passive participants inform the TSO of their deviations from the contracted schedule. Such deviations are priced under a dual-price imbalance settlement scheme, i.e., a different price for positive (extra-production) and negative (under-production) deviations (\cite{morales2013integrating, morales2010short}). If the passive participant generates a positive deviation $q^{(+)}_k$, the associated income is $\lambda^{(+)}_k q^{(+)}_k$, where $\lambda^{(+)}_k$ is given by
    \begin{equation*}
			\lambda^{(+)}_k = \begin{cases} \lambda^{\mathrm{DA}}_k, & \mbox{if} \hspace{.2cm} \lambda^{\mathrm{BA}}_k \geq \lambda^{\mathrm{DA}}_k \\ 
			\lambda^{\mathrm{BA}}_k, & \mbox{otherwise}\end{cases} 
    \end{equation*}
i.e., the least convenient price between $\lambda^{\mathrm{DA}}_k$ and $\lambda^{\mathrm{BA}}_k$. Similarly, a negative deviation $q^{(-)}_k$ generates a revenue $-\lambda^{(-)}_k q^{(-)}_k$, with $\lambda^{(-)}_k$ given by
    \begin{equation*}
    	\lambda^{(-)}_k = \begin{cases} \lambda^{\mathrm{BA}}_k, & \mbox{if} \hspace{.2cm} \lambda^{\mathrm{BA}}_k \geq \lambda^{\mathrm{DA}}_k\\ 
			\lambda^{\mathrm{DA}}_k, & \mbox{otherwise}\end{cases} 
    \end{equation*}
    i.e., again, the least convenient price between $\lambda^{\mathrm{DA}}_k$ and $\lambda^{\mathrm{BA}}_k$. For a more extended discussion on the dual-pricing mechanism, we refer the interested reader to \cite{dent2011opportunity} and \cite{skajaa2015intraday}.
		
	This balancing market structure (i.e., uniform pricing for settling active offers and dual-pricing for passive deviations) is widely used across Europe, e.g., in Spain, Portugal, and Denmark among other countries \cite{wang2015review}. Our novel \textit{Active/Passive} market model thus adapts to one of the major market contexts.
		
	Notice that the contracted schedule of active participants is rigid and deviations are not allowed in our market model. In practice, small output deviations from the schedule are sometimes tolerated, mainly in balancing markets with shorter trading intervals (e.g., 5 minutes). However, this relaxation should not be considered in the offering strategy of an active participant, as it is usually introduced to account for possible control errors in the real-time operation of the power units.

\subsection{VPP structure}\label{subsec:VPP_Structure}
		
\begin{figure}[t!]
	\centering
	\includegraphics[width=0.7\columnwidth]{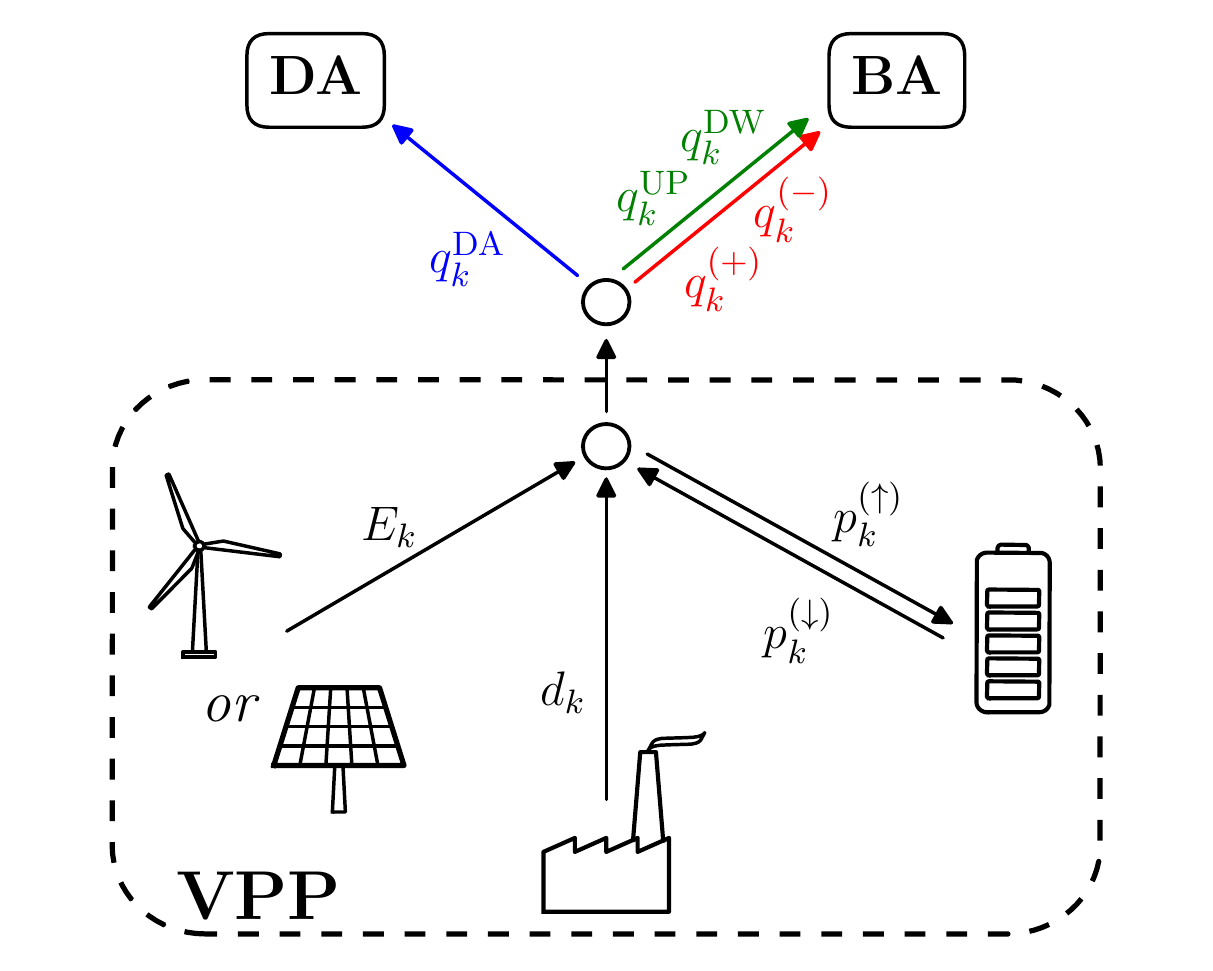}
	\caption{Illustration of the VPP structure highlighting the energy exchanged with the electricity market platform.}             \label{fig:VPP_structure}
\end{figure}

	We consider a VPP composed of a stochastic power unit (either wind or solar), a conventional thermal unit, and an electric energy storage. The structure of the VPP is illustrated in Fig.~\ref{fig:VPP_structure}. The power production of the stochastic unit and the thermal unit are denoted by $E_k$ and $d_k$, respectively. The storage unit produces energy when discharging and consumes energy during the charging phase; the amount of charging and discharging power is denoted by $p^{(\uparrow)}_k$ and $p^{(\downarrow)}_k$, respectively. The total amount of energy production (or consumption) of the VPP has to match the amount of energy exchanged with the electricity market platform. The energy quantity contracted in the day-ahead market is denoted by $q^{\mathrm{DA}}_k$. We then indicate with $q^{\mathrm{UP}}_k$ and $q^{\mathrm{DW}}_k$ the upward and downward adjustments in the balancing market, respectively, which are associated with an active participation at the balancing stage. Alternatively, the VPP can create a positive $q^{(+)}_k$ or negative $q^{(-)}_k$ deviation in the real-time. Thus, under the \textit{Active/Passive} participation model, the ``deterministic'' energy balance at the hourly interval $k$ between the VPP production (or consumption) and the quantity exchanged with the market platform can be expressed by
	\begin{equation}\label{eq:energy_balance_VPP}
		q^{\mathrm{DA}}_k + q^{\mathrm{UP}}_k - q^{\mathrm{DW}}_k + q^{(+)}_k - q^{(-)}_k = E_k + d_k + p^{(\downarrow)}_k - p^{(\uparrow)}_k. \nonumber
	\end{equation}
The VPP is assumed price-taker in both the day-ahead and the balancing market. Accordingly, the market prices within its offering strategy are exogenous and uncertain, and modeled by means of a set of scenarios. We also assume that the different technologies within the VPP are owned by the same company, thus, in this paper we do not investigate how to eventually redistribute the VPP profit among its components.
		
\subsection{Scenario Generation}\label{subsec:Scenario_Generation}
		
	To derive the offering strategy, the VPP is provided with an input set of scenarios describing the evolution of uncertain market prices and power production from the wind or PV power unit. We generate such scenarios starting from probabilistic forecasts. The probabilistic forecasts for the day-ahead and the balancing market prices are obtained through the fundamental market model proposed in \cite{mazzi2017price}, where parametrized supply and demand curves are used to simulate the market clearing mechanism. Then, by introducing uncertainty in one or more parameters of the two curves, we obtain probabilistic forecasts of the day-ahead and balancing market prices. For wind and PV power production we instead use the dataset of probabilistic forecasts, respectively, from \cite{mazzi2017price} and \cite{pierro2016multi}.
		
	Probabilistic forecasts describe an estimate of the random variable density function for each look-ahead time, without any inter-temporal correlation. To include temporal dependencies, starting from the probabilistic forecasts we generate a set of trajectories following the methodology presented in \cite{pinson2009probabilistic} and \cite{pinson2012evaluating}. In brief, series of forecast errors are converted into a multivariate Gaussian random variable, and a unique covariance matrix is used to describe the interdependence structure. We model this covariance matrix using an exponential covariance function \cite{pinson2012evaluating} where the exponential parameter controls the correlation among the lead times. Accordingly, the day-ahead market prices $\lambda^{\mathrm{DA}}_k$ are represented by the set of trajectories $\{ \lambda^{\mathrm{DA}}_{ik} : i \in I, k \in K \}$. Then, for each day-ahead scenario $i$, the balancing market prices $\lambda^{\mathrm{BA}}_k$ are modeled using a set $J$ of scenarios $\{ \lambda^{\mathrm{BA}}_{ijk} : i \in I, j \in J, k \in K \}$. Finally, the uncertain power production $E_k$ from the stochastic unit is represented by the set of trajectories $\{ E_{\omega k} : \omega \in W, k \in K \}$.
		
	The number of scenarios needed to accurately represent continuous random variables or stochastic processes is usually large, leading to intractable stochastic programs. Therefore, we use the technique of \cite{growe2003scenario} to reduce the number of scenarios while preserving most of the stochastic information.

\section{Optimal Offering Strategy through Multi-Stage Stochastic Programming}\label{sec:Optimal_Offering_Strategy}
\begin{figure}[b!]
	\centering
	\includegraphics[width=0.84\columnwidth]{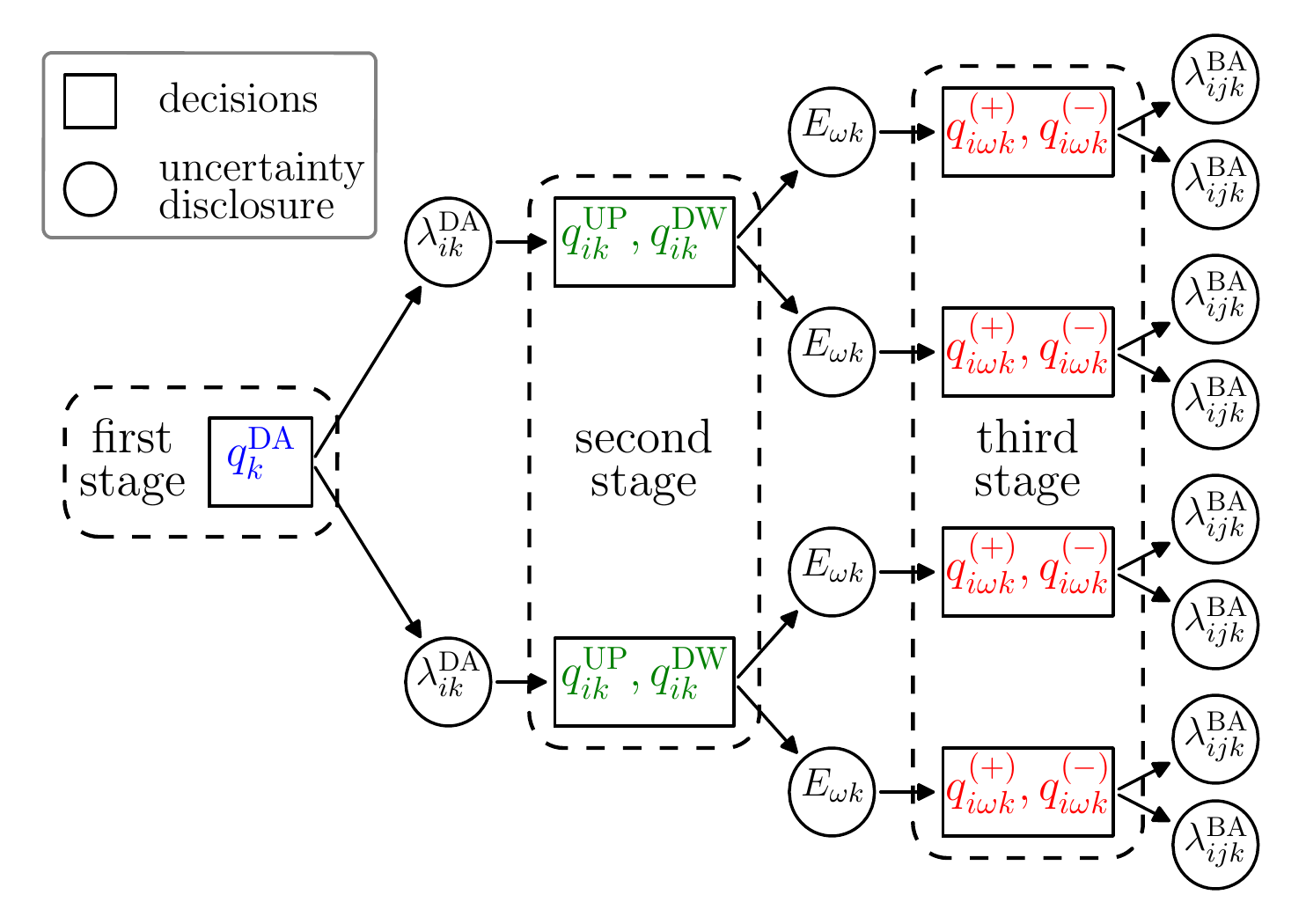}
	\caption{Representation of the multi-stage stochastic programming setup.} \label{fig:Stoch_Prog_FMwork}
\end{figure}		
	At noon, the VPP submits its day-ahead market offers to the market operator, aiming to maximize its total expected profit. While determining the optimal day-ahead market offers, the VPP also takes into account the uncertainty, and endogenously models the future decisions in the balancing market. This results in a three-stage stochastic programming framework that we illustrate in Fig.~\ref{fig:Stoch_Prog_FMwork}. The day-ahead quantity offers, $q^{\mathrm{DA}}_k$, are modeled as first-stage decisions. Then, at the second stage, the VPP owner makes the decision of being active or passive. Specifically, at the second stage, the day-ahead market price $\lambda^{\textrm{DA}}_{ik}$ for hour $k$ is known and the VPP decides whether it will be active or passive in hour $k$. If the chosen behavior is active, then, still at the second stage, the VPP submits its offering curves to the balancing market for up- or down-regulation. Therefore, the up- and down-regulation adjustments $q^{\mathrm{UP}}_{ik}$ and $q^{\mathrm{DW}}_{ik}$ are second-stage decision variables that are chosen after the day-ahead market prices $\lambda^{\mathrm{DA}}_k$ realize but before the power production $E_{k}$ is known. In contrast, if the chosen behavior is passive, then the VPP does not offer regulating energy to the balancing market but instead it will settle the imbalances in the third stage, which means that the positive and negative deviations $q^{(+)}_{i\omega k}$ and $q^{(-)}_{i\omega k}$ are third-stage decisions that are made after the disclosure of $\lambda^{\mathrm{DA}}_k$ and $E_{k}$.
		
	Following the methodology presented in \cite{conejo2010decision} and \cite{mazzi2017price}, we make the day-ahead quantities $q^{\mathrm{DA}}_k$ scenario $i$ dependent (i.e., $q^{\mathrm{DA}}_k \rightarrow q^{\mathrm{DA}}_{ik}$) to build the non-decreasing offer curves. Despite being built using scenario-dependent price-quantity offers, the curves adapt to any realization of the uncertainty and are, in fact, scenario-independent (the non-anticipativity structure of our stochastic program is not violated). Similarly, the up- and down-regulation adjustments $q^{\mathrm{UP}}_{ik}$ and $q^{\mathrm{DW}}_{ik}$ are made index $j$ dependent (i.e., $q^{\mathrm{UP}}_{ik} \rightarrow q^{\mathrm{UP}}_{ijk}$ and $q^{\mathrm{DW}}_{ik} \rightarrow q^{\mathrm{DW}}_{ijk}$) to derive the offer curves in the balancing market. Finally, we make the VPP operational variables dependent on indices $i$, $j$, and $\omega$ (e.g., $d_k \rightarrow d_{ij\omega k}$ and $p^{(\uparrow)}_k \rightarrow p^{(\uparrow)}_{ij \omega k}$). To derive its day-ahead market offers, the VPP solves the following optimization model
\begin{subequations}
	\label{eq:Off_Strategy}
	\begin{align}
	\max_{\Xi}
	& \quad \sum_{k} \hat{\rho}^{\mathrm{DA}}_k + \hat{\rho}^{\mathrm{Act}}_k + \hat{\rho}^{\mathrm{Pas}}_k - \hat{c}_k \label{eq:Off_Strategy_Obj}\\
	\begin{split} 
	\text{s.t.} 
	& \quad q^{\textrm{DA}}_{ik} + q^{\textrm{UP}}_{ijk} - q^{\textrm{DW}}_{ijk} + q^{(+)}_{i\omega k} - q^{(-)}_{i\omega k} = \\ 
			& \qquad E_{\omega k} + d_{ij \omega k} + p^{(\uparrow)}_{ij \omega k} - p^{(\downarrow)}_{ij \omega k} , \,\,\, \forall i, \forall j, \forall \omega, \forall k \label{eq:Off_Strategy_Balance}  \end{split}\\ 
			& \quad \left( \hat{\rho}^{\mathrm{DA}}_k, q^{\textrm{DA}}_{ik} \right) \in \Pi^{\mathrm{DA}}, \quad \forall i, \forall k \label{eq:Off_Strategy_Pi_DA}\\
			& \quad \left( \hat{\rho}^{\mathrm{Act}}_k, q^{\textrm{UP}}_{ijk}, q^{\textrm{DW}}_{ijk} \right) \in \Pi^{\mathrm{Act}}, \quad \forall i, \forall j, \forall k \label{eq:Off_Strategy_Pi_Act}\\
			& \quad \left( \hat{\rho}^{\mathrm{Pas}}_k, q^{(+)}_{ijk}, q^{(-)}_{ijk} \right) \in \Pi^{\mathrm{Pas}}, \quad \forall i, \forall \omega, \forall k \label{eq:Off_Strategy_Pi_Pas}\\
			& \quad \left( d_{ij\omega k}, p^{(\uparrow)}_{ij\omega k}, p^{(\downarrow)}_{ij\omega k} \right) \in \Omega, \quad \forall i, \forall j, \forall \omega, \forall k \label{eq:Off_Strategy_Omega}\\
			& \quad  \hat{c}_k = h \hspace{-1.8pt} \left( \left\{ d_{ij\omega k}, \hspace{3pt} \forall i, \forall j, \forall \omega \right\} \right) , \quad \forall k \label{eq:Off_Strategy_g}\\
			& \quad \left( q^{\textrm{UP}}_{ijk}, q^{\textrm{DW}}_{ijk},q^{(+)}_{i\omega k} ,q^{(-)}_{i\omega k}\right) \in \Gamma,\quad  \forall i, \forall j, \forall \omega, \forall k \, \label{eq:Off_Strategy_Gamma}
	\end{align}
	\end{subequations}\\
	\begin{multline*}
		\text{where}\quad \Xi = \big\{ \hat{\rho}^{\mathrm{DA}}_k, \hat{\rho}^{\mathrm{Act}}_k, \hat{\rho}^{\mathrm{Pas}}_k, \hat{c}_k,  q^{\textrm{DA}}_{ik} , q^{\textrm{UP}}_{ijk} , q^{\textrm{DW}}_{ijk} , \\  q^{(+)}_{i\omega k} , q^{(-)}_{i\omega k}, d_{ij \omega k} , p^{(\uparrow)}_{ij \omega k} , p^{(\downarrow)}_{ij \omega k} \big\}.
	\end{multline*}
		
	The objective function \eqref{eq:Off_Strategy_Obj} maximizes the VPP expected revenues considering both the day-ahead and balancing markets. Constraint \eqref{eq:Off_Strategy_Balance} imposes the energy balance between the VPP production (or consumption) and the energy exchanged with the electricity market. The sets of constraints \eqref{eq:Off_Strategy_Pi_DA}-\eqref{eq:Off_Strategy_Gamma} are expressed and discussed in detail in the following.

\subsection{Linear Formulation of $ \Pi^{\mathrm{DA}} $}\label{subsec:Linear_Expression_Pi_DA}
		
	The set of constraints \eqref{eq:Off_Strategy_Pi_DA}, denoted by $\Pi^{\mathrm{DA}}$, computes the expected profit from the day-ahead market $\hat{\rho}^{\mathrm{DA}}_k$ and includes constraints on the day-ahead offer curve. It is written as
	\begin{subequations}
	\label{eq:Pi_DA_Linear}
	\begin{align}
			& \hat{\rho}^{\mathrm{DA}}_k = \sum_i \pi^{\mathrm{DA}}_i \lambda^{\mathrm{DA}}_{ik} q^{\mathrm{DA}}_{ik}, \quad \forall k \label{eq:Pi_DA_Linear_rhoDA}\\ 
			& q^{\mathrm{DA}}_{ik} \geq q^{\mathrm{DA}}_{i' k} \quad \mbox{if} \quad \lambda^{\mathrm{DA}}_{ik} \geq \lambda^{\mathrm{DA}}_{i' k}, \quad \forall i, \forall i', \forall k \label{eq:Pi_DA_Linear_nondec}\\ 
			& q^{\mathrm{DA}}_{ik} = q^{\mathrm{DA}}_{i' k} \quad \mbox{if} \quad \lambda^{\mathrm{DA}}_{ik} = \lambda^{\mathrm{DA}}_{i' k}, \quad \forall i, \forall i', \forall k \label{eq:Pi_DA_Linear_nonant}\\ 
			& - \overline{P}^{(\uparrow)} \leq q^{\mathrm{DA}}_{ik} \leq \overline{D} + \overline{E} +  \overline{P}^{(\downarrow)}, \quad \forall i, \forall i', \forall k \label{eq:Pi_DA_Linear_qDA}.
	\end{align}
	\end{subequations}
Constraint \eqref{eq:Pi_DA_Linear_rhoDA} yields the expected income of the VPP associated with the day-ahead market offer curves. Constraints \eqref{eq:Pi_DA_Linear_nondec} and \eqref{eq:Pi_DA_Linear_nonant} force the offer curves to be, respectively, non-decreasing and non-anticipative. Finally, constraints \eqref{eq:Pi_DA_Linear_qDA} restrict the day-ahead offer quantities to the VPP capacity.

\subsection{Linear Formulation of $ \Pi^{\mathrm{Act}} $}\label{subsec:Linear_Expression_Pi_Act}
		
	The constraints \eqref{eq:Off_Strategy_Pi_Act}, denoted by $\Pi^{\mathrm{Act}}$, yield the expected profit $\hat{\rho}^{\mathrm{Act}}_k$ from an active participation in the balancing market, and comprises constraints on the offer curves in the balancing market. They are formulated as
\begin{subequations}
	\label{eq:Pi_Act_Linear}
	\begin{align}
			& \hat{\rho}^{\mathrm{Act}}_k = \sum_{ij} \pi^{\mathrm{DA}}_i \pi^{\mathrm{BA}}_j \lambda^{\mathrm{BA}}_{ijk} \left( q^{\mathrm{UP}}_{ijk} - q^{\mathrm{DW}}_{ijk} \right)  , \quad \forall k \label{eq:Pi_Act_Linear_rho}\\ 
			& q^{\mathrm{UP}}_{ijk} \geq q^{\mathrm{UP}}_{i j' k} \quad \mbox{if} \quad \lambda^{\mathrm{BA}}_{ijk} \geq \lambda^{\mathrm{BA}}_{i j' k}, \quad \forall i, \forall j,  \forall j', \forall k \label{eq:Pi_Act_Linear_nondecUP}\\ 
			& q^{\mathrm{UP}}_{ijk} = q^{\mathrm{UP}}_{i j' k} \quad \mbox{if} \quad \lambda^{\mathrm{BA}}_{ijk} = \lambda^{\mathrm{BA}}_{i j' k}, \quad \forall i, \forall j, \forall j', \forall k \label{eq:Pi_Act_Linear_nonantUP}\\ 
			& q^{\mathrm{DW}}_{ijk} \leq q^{\mathrm{DW}}_{i j' k} \quad \mbox{if} \quad \lambda^{\mathrm{BA}}_{ijk} \geq \lambda^{\mathrm{BA}}_{i j' k}, \quad \forall i, \forall j,  \forall j', \forall k \label{eq:Pi_Act_Linear_nondecDW}\\ 
			& q^{\mathrm{DW}}_{ijk} = q^{\mathrm{DW}}_{i j' k} \quad \mbox{if} \quad \lambda^{\mathrm{BA}}_{ijk} = \lambda^{\mathrm{BA}}_{i j' k}, \quad \forall i, \forall j, \forall j', \forall k \label{eq:Pi_Act_Linear_nonantDW}\\ 
			& q^{\mathrm{UP}}_{ijk} = 0 \;\; \mbox{if down-regulation in ($i,j$) at $k$} , \quad \forall i, \forall j, \forall k \label{eq:Pi_Act_Linear_qUP0}\\ 
			& q^{\mathrm{DW}}_{ijk} = 0 \;\; \mbox{if up-regulation in ($i,j$) at $k$}, \quad \forall i, \forall j, \forall k \label{eq:Pi_Act_Linear_qDW0}\\ 
			& q^{\mathrm{UP}}_{ijk}, q^{\mathrm{DW}}_{ijk} \geq 0, \quad \forall i, \forall j, \forall k. \label{eq:Pi_Act_Linear_qUPDWpos}
	\end{align}
\end{subequations}
	Constraint \eqref{eq:Pi_Act_Linear_rho} evaluates the expected revenue from the submission of offer curves in the balancing market as an active participant. Constraints \eqref{eq:Pi_Act_Linear_nondecUP} and \eqref{eq:Pi_Act_Linear_nonantUP} ensure, respectively non-decreasing shape and non-anticipativity of the up-regulation offer curve. Similarly, constraints \eqref{eq:Pi_Act_Linear_nondecDW} and \eqref{eq:Pi_Act_Linear_nonantDW} do the same for the down-regulation offer curve. 
    Constraints \eqref{eq:Pi_Act_Linear_qUP0} and \eqref{eq:Pi_Act_Linear_qDW0} impose the offer type (up- or down-regulation) based on the direction of the system imbalance, information that is embedded in each scenario $(i,j)$ at interval $k$. Nevertheless, in the case study of Section V, we make the simplifying assumption (also made in [25]) that the TSO requires down-regulation in ($i,j$) at $k$ if and only if $\lambda^{\mathrm{BA}}_{ijk} < \lambda^{\mathrm{DA}}_{ik}$, and up-regulation if and only if $\lambda^{\mathrm{BA}}_{ijk} > \lambda^{\mathrm{DA}}_{ik}$, that is, we assume that the system's need for up- or down-regulation can be directly inferred from the difference between $\lambda^{\mathrm{DA}}_{ik}$ and $\lambda^{\mathrm{BA}}_{ijk}$. Note that situations in which this is not true are infrequent in real-word dual-price balancing markets. Furthermore, the benefits the VPP could make from trying to anticipate these rare events when offering in the day-ahead market would not compensate for the high costs such an endeavor would involve in terms of computational and modeling effort and the risk of a wrong guess.
	Finally, constraint \eqref{eq:Pi_Act_Linear_qUPDWpos} enforces $q^{\mathrm{UP}}_{ijk}$ and $q^{\mathrm{DW}}_{ijk}$ to be non-negative variables.

\subsection{Linear Formulation of $ \Pi^{\mathrm{Pas}} $}\label{subsec:Linear_Expression_Pi_Pas}
		
	The constraints \eqref{eq:Off_Strategy_Pi_Pas}, denoted by $\Pi^{\mathrm{Pas}}$, give the expected profit $\hat{\rho}^{\mathrm{Pas}}_k$ associated with a passive participation in the balancing market. They are formulated as
\begin{subequations}
	\label{eq:Pi_Pas_Linear}
	\begin{align}
			& \hat{\rho}^{\mathrm{Pas}}_k = \sum_{ij\omega } \pi^{\mathrm{DA}}_i \pi^{\mathrm{BA}}_j \pi^{E}_{\omega} \left( \lambda^{\mathrm{(+)}}_{ijk} q^{(+)}_{i\omega k} - \lambda^{\mathrm{(-)}}_{ijk} q^{(-)}_{i\omega k} \right)  , \quad \forall k \label{eq:Pi_Pas_Linear_rho}\\ 
			& q^{(+)}_{i\omega k}, q^{(-)}_{i\omega k} \geq 0, \quad \forall i, \forall \omega, \forall k \label{eq:Pi_Pas_Linear_qPMpos}
	\end{align}
\end{subequations}
where $\lambda^{\mathrm{(+)}}_{ijk} = \min \hspace{-1pt} \big( \lambda^{\mathrm{BA}}_{ijk},\lambda^{\mathrm{DA}}_{ik} \big) $ and $\lambda^{\mathrm{(-)}}_{ijk} = \max \hspace{-1pt} \left( \lambda^{\mathrm{BA}}_{ijk},\lambda^{\mathrm{DA}}_{ik} \right) $, in accordance with the dual-price imbalance settlement scheme \cite{morales2013integrating, morales2010short}. Constraint \eqref{eq:Pi_Pas_Linear_rho} computes the expected income from a passive participation in the balancing stage and accounts for the imbalances created. Constraint \eqref{eq:Pi_Pas_Linear_qPMpos} ensures that $q^{(+)}_{i\omega k} $ and $ q^{(-)}_{i\omega k}$ are non-negative variables.

\subsection{MILP Formulation of $ \Omega $}\label{subsec:Linear_Expression_Omega}
		
	The constraints \eqref{eq:Off_Strategy_Omega}, denoted by $\Omega$, establish the feasible operating region of the VPP and are formulated as
\begin{subequations}
	\label{eq:Omega_Linear}
	\begin{align}
	& \ell_{ij \omega k} = \ell_{ij \omega (k-1)} + \eta \hspace{1pt} p^{(\uparrow)}_{ij\omega k} - p^{(\downarrow)}_{ij\omega k}   , \quad \forall i, \forall j, \forall \omega, \forall k \label{eq:Omega_Linear_lev}\\ 
	&  \underline{L} \leq \ell_{ij \omega k} \leq \overline{L}   , \quad \forall i, \forall j, \forall \omega, \forall k \label{eq:Omega_Linear_levMn}\\ 
	&  0 \leq p^{(\uparrow)}_{ij \omega k} \leq \overline{P}^{(\uparrow)}   , \quad \forall i, \forall j, \forall \omega, \forall k \label{eq:Omega_Linear_pchMn}\\ 
	&  0 \leq p^{(\downarrow)}_{ij \omega k} \leq \overline{P}^{(\downarrow)}   , \quad \forall i, \forall j, \forall \omega, \forall k \label{eq:Omega_Linear_pdsMn}\\ 
	&  u_{ij \omega k} \hspace{1pt} \underline{D} \leq d_{ij \omega k} \leq u_{ij \omega k} \hspace{1pt} \overline{D}   , \quad \forall i, \forall j, \forall \omega, \forall k \label{eq:Omega_Linear_dMn}\\ 
	&  d_{ij \omega k} - d_{ij \omega (k-1)} \leq R^{\mathrm{UP}}   , \quad \forall i, \forall j, \forall \omega, \forall k \label{eq:Omega_Linear_drampUP}\\ 
	&  d_{ij \omega (k-1)} - d_{ij \omega k} \leq R^{\mathrm{DW}}   , \quad \forall i, \forall j, \forall \omega, \forall k \label{eq:Omega_Linear_drampDW}\\ 
	& u_{ij \omega k} \in \{0,1\}, \quad \forall i, \forall j, \forall \omega, \forall k \label{eq:Omega_Linear_ubin}.
	\end{align}
\end{subequations}
Constraint \eqref{eq:Omega_Linear_lev} represents the energy balance of the storage unit. Constraint \eqref{eq:Omega_Linear_levMn} forces the level of energy in the storage $\ell_{ij\omega k}$ to lie between its minimum and maximum limits. Similarly, constraints \eqref{eq:Omega_Linear_pchMn} and \eqref{eq:Omega_Linear_pdsMn} do the same for $p^{(\uparrow)}_{ij\omega k}$ and $p^{(\downarrow)}_{ij\omega k}$, respectively. Constraint \eqref{eq:Omega_Linear_dMn} imposes the thermal unit to operate within its minimum output and its capacity when on-line (i.e., $u_{ij\omega k}=1$) and not to produce when off-line (i.e., $u_{ij\omega k}=0$). Constraints \eqref{eq:Omega_Linear_drampUP} and \eqref{eq:Omega_Linear_drampDW} enforce, respectively, the upward and downward ramping limitations of the thermal unit. Finally, constraint \eqref{eq:Omega_Linear_ubin} sets the commitment status $u_{ij\omega k}$ of the thermal unit as a binary variable. Constraint \eqref{eq:Omega_Linear_lev} requires the initial level of the storage as input, and constraints \eqref{eq:Omega_Linear_drampUP}-\eqref{eq:Omega_Linear_drampDW} need the initial production level of the thermal unit.
		
	Richer models for the feasible operating region of a dispatchable unit exist in the literature. However, to keep the focus on the \textit{Active/Passive} offering strategy and its formulation intuitive, we chose to capture the main operating constraints of the unit limiting the level of details of the model.

\subsection{MILP Formulation of $ h(\cdot) $}\label{subsec:Linear_Expression_g}
		
	Constraint \eqref{eq:Off_Strategy_g} computes the expected production cost $\hat{c}_{k}$ associated with the thermal unit. A possible mixed-integer linear programming formulation of this cost function is \eqref{eq:h_Linear}
	\begin{equation}\label{eq:h_Linear}
		\hat{c}_k = \sum_{ij\omega } \pi^{\mathrm{DA}}_i \pi^{\mathrm{BA}}_j \pi^{E}_{\omega} \left( C_0 u_{ij\omega k} + C \hspace{1pt} d_{ij\omega k} \right), \,\,\forall k 
	\end{equation}
where $C_0$ is the fixed cost incurred when the unit is on, and $C$ is the marginal production cost of the unit.
		
\subsection{MILP Formulation of $ \Gamma $}\label{subsec:Linear_Expression_Gamma}
		
	The set of constraints \eqref{eq:Off_Strategy_Gamma}, denoted by $\Gamma$, enforces complementarity between the active and passive participation in the balancing market. $\Gamma$ can be formulated as
	\begin{subequations}
		\label{eq:Gamma_Linear}
		\begin{align}
			& q^{\mathrm{UP}}_{ijk} + q^{\mathrm{DW}}_{ijk} \leq \epsilon_{ik} \overline{M}  , \quad \forall i, \forall j, \forall k \label{eq:Gamma_Linear_UPDW}\\ 
			& q^{(+)}_{i\omega k} + q^{(-)}_{i\omega k} \leq \left( 1-\epsilon_{ik} \right) \overline{M}  , \quad \forall i, \forall \omega, \forall k \label{eq:Gamma_Linear_PM}\\ 
			& \epsilon_{i k} \in \{0,1\}, \quad \forall i, \forall k \label{eq:Gamma_Linear_epsbin}.
		\end{align}
	\end{subequations}
Constraints \eqref{eq:Gamma_Linear_UPDW} and \eqref{eq:Gamma_Linear_PM} force the VPP to be in only one state between active (i.e., $\epsilon_{ik}=1$) and passive (i.e., $\epsilon_{ik}=0$) in the balancing market, through the so-called \textit{big-M} approach. A natural and sensible choice for the parameter $\overline{M}$ can be
\begin{equation*}
	\overline{M} := \overline{E} + \overline{D} + \overline{P}^{(\uparrow)} + \overline{P}^{(\downarrow)}.
\end{equation*}
Finally, constraint \eqref{eq:Gamma_Linear_epsbin} forces the variables $\epsilon_{ik}$ to be binary.
		
\section{Example}\label{sec:Example}
        
	In this section, we present a stylized example to better illustrate the potential contribution of enabling \textit{Active/Passive} offers in the balancing market. We analyze the offering strategy of a VPP composed of a wind farm of capacity $\overline{E}=40$ MW and a conventional generation unit of capacity $\overline{D}=25$ MW, minimum output $\underline{D}=0$ MW, and marginal cost $C=31$ \euro/MWh ($C_0=0$ \euro). This example considers 2 trading intervals ($k_1$,$k_2$), 1 day-ahead price scenario, 2 balancing-market price scenarios ($j_1$,$j_2$), and 2 wind power production scenarios ($\omega_1$,$\omega_2$). The values of $\lambda^{\mathrm{DA}}_k$, $\lambda^{\mathrm{BA}}_{jk}$ (and the associated $\lambda^{(+)}_{jk}$ and $\lambda^{(-)}_{jk}$), and $E_{\omega k}$ are given in Table \ref{tab:PW_Scen_Example}. Note that the balancing scenarios $j_1$ and $j_2$ have the same probability of realization, likewise the wind production scenarios $\omega_1$ and $\omega_2$, i.e., $\pi^{\mathrm{BA}}_{j_1}=\pi^{\mathrm{BA}}_{j_2}=0.5$ and $\pi^{\mathrm{E}}_{\omega_1}=\pi^{\mathrm{E}}_{\omega_2}=0.5$.

\begin{table}[h!]
	\centering
	\caption{Prices and wind production scenarios.}
	\label{tab:PW_Scen_Example}
	\begin{tabular}{cccccccccc}
		\toprule
		& $\lambda^{\mathrm{DA}}_k$ & \multicolumn{2}{c}{$\lambda^{\mathrm{BA}}_{jk}$} & \multicolumn{2}{c}{$\lambda^{(+)}_{\omega k}$} & \multicolumn{2}{c}{$\lambda^{(-)}_{\omega k}$} & \multicolumn{2}{c}{$E_{\omega k}$} \\
		& (\euro/MWh)               & \multicolumn{2}{c}{(\euro/MWh)}                  & \multicolumn{2}{c}{(\euro/MWh)}                & \multicolumn{2}{c}{(\euro/MWh)}                & \multicolumn{2}{c}{(MWh)}          \\
		&                           & $j_1$                   & $j_2$                  & $j_1$                  & $j_2$                 & $j_1$                  & $j_2$                 & $\omega_1$       & $\omega_2$      \\ \midrule
		$k_1$ & 25                        & 26                      & 23                     & 25                     & 23                    & 26                     & 25                    & 5                & 18              \\
		$k_2$ & 29                        & 19                      & 37                     & 19                     & 29                    & 29                     & 37                    & 9                & 15              \\
     \bottomrule
	\end{tabular}
\end{table}
			
	The optimal VPP offers under the current balancing market framework, i.e.\ passive-only, are shown in Table \ref{tab:Decisions_Pas}. 
    \begin{table}[b!]
	\centering
	\caption{Decisions with Passive strategy.}
	\label{tab:Decisions_Pas}
	\begin{tabular}{ccccccccccc}
	\toprule
	& $q^{\mathrm{DA}}_k$ & $\epsilon_{k}$ & \multicolumn{2}{c}{$q^{\mathrm{UP}}_{jk}$} & \multicolumn{2}{c}{$q^{\mathrm{DW}}_{jk}$} & \multicolumn{2}{c}{$q^{(+)}_{\omega k}$} & \multicolumn{2}{c}{$q^{(-)}_{\omega k}$} \\
	& (MWh)               &                & \multicolumn{2}{c}{(MWh)}                  & \multicolumn{2}{c}{(MWh)}                  & \multicolumn{2}{c}{(MWh)}                & \multicolumn{2}{c}{(MWh)}                \\
	&                     &                & $j_1$                & $j_2$               & $j_1$                & $j_2$               & $\omega_1$          & $\omega_2$         & $\omega_1$          & $\omega_2$         \\ \midrule
	$k_1$ & 18                  & 0              & 0                    & 0                   & 0                    & 0                   & 0                   & 0                  & 13                  & 0                  \\
	$k_2$ & 15                  & 0              & 0                    & 0                   & 0                    & 0                   & 0                   & 0                  & 0                   & 0                  \\
                        \bottomrule
	\end{tabular}
\end{table}
    The VPP deviates -13 MWh at $k_1$ if scenario $\omega_1$ realizes. Differently, at interval $k_2$ it relies on the dispatchable units to balance the wind production uncertainty. Note that this behavior is linked to how different $\lambda^{\mathrm{BA}}_{jk}$ is expected to be from $\lambda^{\mathrm{DA}}_{k}$. If these are ``close'' as in $k_1$, it may be convenient to use the market to compensate for the uncertainty in $E_{\omega k}$. At $k_1$, $E_{\omega k}$ will be either 5 or 18 MWh, and the VPP schedules $q^{\mathrm{DA}}_{k_1}=18$ MWh, receiving a day-ahead income $\hat{\rho}^{\mathrm{DA}}_{k_1}$ of 450 \euro\ ($\lambda^{\mathrm{DA}}_{k_1} \times q^{\mathrm{DA}}_{k_1}$). In scenario $\omega_2$ there is no need for balancing, while in $\omega_1$ the VPP is short of 13 MWh. In this case, the VPP can either deviate -13 MWh or generate 13 MWh with the conventional unit. A deviation of -13 MWh is priced at $26$ \euro/MWh under scenario $j_1$ or $25$ \euro/MWh under scenario $j_2$, with an associated income of $-338$ \euro$\,$and $-325$ \euro, respectively. With the associated probabilities this results in $\hat{\rho}^{\textrm{Pas}}_{k_1}$ of -165.75 \euro\ ($\pi^{\mathrm{E}}_{\omega_1}\times(\pi^{\mathrm{BA}}_{j_1}\lambda^{(-)}_{j_1 k_1}+\pi^{\mathrm{BA}}_{j_2}\lambda^{(-)}_{j_2 k_1})\times q^{(-)}_{\omega_1 k_1}$), yielding a total expected profit of 284.25 \euro. Using the thermal unit to produce these 13 MWh costs 403 \euro, yielding an expected profit of 248.25 \euro, which is less convenient. At $k_2$, the balancing price $\lambda^{\mathrm{BA}}_{jk}$ is likely to be ``far'' from $\lambda^{\mathrm{DA}}_{k}$, and the VPP contracts $q^{\mathrm{DA}}_{k_2}=15$ MWh, which gives a day-ahead income of $435$ \euro. When $\omega_2$ realizes the wind farm produces 15 MWh, thus matching the day-ahead quantity offer. Under scenario $\omega_1$, $E_{\omega k}$ is 9 MWh and the VPP can either produce 6 MWh with the conventional unit, or deviate in the balancing market. The expected cost of producing 6 MWh is 93 \euro. 
Differently, relying on the balancing market yields an expected income $\hat{\rho}^{\textrm{Pas}}_{k_3}$ of -111 \euro, associated with $q^{(-)}_{\omega_1  k_2}=6$ MWh. Then, the decision of producing using the thermal unit is more convenient. Note that if $\lambda^{\mathrm{BA}}_{j_2 k_2}$ would have been ``closer'' to $\lambda^{\mathrm{DA}}_{k_2}$, e.g., $\lambda^{\mathrm{BA}}_{j_2 k_2} = 30$ \euro/MWh, then it would have been more profitable to settle $q^{(-)}_{\omega_1  k_2}=6$ MWh in the balancing market.

	Under the \textit{Active/Passive} model (see Table \ref{tab:Decisions_ActPas}), the VPP offers regulation to the market during interval $k_2$, i.e., when it does not rely on the market to balance its position. At $k_2$, the VPP increases its day-ahead offer in 19 MWh (with respect to what the VPP does in passive-only mode), which gives an additional day-ahead income of 551 \euro. Producing such 19 MWh is needed only in scenario $j_2$ and the associated expected cost is 294.5 \euro. Differently, under scenario $j_1$ it submits a down-regulation offer $q^{\mathrm{DW}}_{j_1 k_2}=19$ MWh, which is priced at 19 \euro/MWh. This generates an expected income $\hat{\rho}^{\mathrm{Act}}_{j_1 k_2}$ of -180.5 \euro\ (negative as it is a down-regulation offer), but avoids the cost of producing the 19 MWh (under scenario $j_1$). In total, this gives an extra expected profit of 76 \euro, i.e.,  $551-294.5-180.5$, at $k_2$ with respect to the passive-only strategy.

\begin{table}[htpb!]
	\centering
	\caption{Decisions with Active/Passive strategy.}
	\label{tab:Decisions_ActPas}
	\begin{tabular}{ccccccccccc}
	\toprule
	& $q^{\mathrm{DA}}_k$ & $\epsilon_{k}$ & \multicolumn{2}{c}{$q^{\mathrm{UP}}_{jk}$} & \multicolumn{2}{c}{$q^{\mathrm{DW}}_{jk}$} & \multicolumn{2}{c}{$q^{(+)}_{\omega k}$} & \multicolumn{2}{c}{$q^{(-)}_{\omega k}$} \\
	& (MWh)               &                & \multicolumn{2}{c}{(MWh)}                  & \multicolumn{2}{c}{(MWh)}                  & \multicolumn{2}{c}{(MWh)}                & \multicolumn{2}{c}{(MWh)}                \\
	&                     &                & $j_1$                & $j_2$               & $j_1$                & $j_2$               & $\omega_1$          & $\omega_2$         & $\omega_1$          & $\omega_2$         \\ \midrule
	$k_1$ & 18                  & 0              & 0                    & 0                   & 0                    & 0                   & 0                   & 0                  & 13                  & 0                  \\
	$k_2$ & 34                  & 1              & 0                    & 0                   & 19                   & 0                   & 0                   & 0                  & 0                   & 0                  \\
    \bottomrule
	\end{tabular}
\end{table}
    
	This small examples shows that both the VPP and the TSO can benefit from enabling the proposed \textit{Active/Passive} strategy in the balancing market. Indeed, the VPP can increase its expected revenues and the TSO can have more regulating flexibility. This example also gives some additional insights. Generally speaking, a balancing market price $\lambda^{\mathrm{BA}}_{jk}$ which is much higher than the day-ahead price $\lambda^{\mathrm{DA}}_{k}$ is typically associated with a significant lack of production in the system, and a $\lambda^{\mathrm{BA}}_{jk}$ much lower than $\lambda^{\mathrm{DA}}_{k}$ to a large excess of production at time $k$. Conversely, when $\lambda^{\mathrm{BA}}_{jk}$ is ``close'' to $\lambda^{\mathrm{DA}}_{k}$, the system imbalance is likely to be small. We have noticed that the VPP is willing to be passive when $\lambda^{\mathrm{BA}}_{jk}$ and $\lambda^{\mathrm{DA}}_{k}$ are ``close'', and active when they are ``far''. Accordingly, the VPP would create deviations when the system imbalance is easier to restore, and provide regulating flexibility when more needed by the TSO.
					
\section{Case Study}\label{sec:Case_Study}
					
	Next we present a case study to test the offering strategy of Section \ref{sec:Optimal_Offering_Strategy}. The aim is to analyze whether the proposed \textit{Active/Passive} balancing participation setup may drive the VPP to offer its flexibility when available.
					
	The scenarios provided as input to the offering model are generated as described in Section \ref{subsec:Scenario_Generation}. First, we generate 300 scenarios for the day-ahead market price $\lambda^{\mathrm{DA}}_{ik}$ and keep the ten most representative ones. Then, for each day-ahead scenario $i$, we randomly sample 300 scenarios of the balancing market price $\lambda^{\mathrm{BA}}_{ik}$ and select the six most significant. Finally, we generate 300 trajectories of the renewable energy production $E_{\omega k}$ (wind or solar power) and keep the five most representative. This results in a scenario tree with 300 branches ($10 \times 6 \times 5$). The parameters of the thermal unit are shown in Table \ref{tab:Thermal_Unit_Data}. Similarly, the characteristics of the storage unit are presented in Table \ref{tab:Storage_Data}. The programs are modeled in \textsc{Python} environment and solved to optimality with \textsc{Gurobi}.
					
	The \textit{Active/Passive} offering strategy is compared against two benchmarks: a \textit{Passive} and an \textit{Active} offering strategy. Based on the \textit{Passive} approach, the VPP is assumed to be always a passive participant in the balancing market. Differently, under the \textit{Active} strategy, the VPP is an active actor in the balancing stage for the entire trading horizon. These two alternative models can be derived from the optimization model \eqref{eq:Off_Strategy} by fixing the binary variables $\epsilon_{ik}=0, \; \forall i, \forall k$ for the \textit{Passive} strategy, or $\epsilon_{ik}=1, \; \forall i, \forall k$ for the \textit{Active} one.
                    
\begin{table}[h!]
	\centering
	\caption{Parameters of the thermal unit.}
	\label{tab:Thermal_Unit_Data}
	\begin{tabular}{cccccc}
	\toprule
	$ \underline{D} $ & $ \overline{D} $ & $ R^{\mathrm{UP}} $ & $ R^{\mathrm{DW}} $ & $ C_0 $ & $ C $         \\
	(MW)     &       (MW)      &        (MW/h)        &        (MW/h)        & (\euro) & (\euro/MWh)   \\ 
	\midrule
	0       &         70       &         30          &          30         & 0     & 45            \\ 
	\bottomrule
	\end{tabular}
\end{table}
    
\begin{table}[h!]
	\centering
	\caption{Parameters of the electric storage unit.}
	\label{tab:Storage_Data}
	\begin{tabular}{cccccc}
	\toprule
	$ \underline{L} $ & $ \overline{L} $ & $ \overline{P}^{(\uparrow)} $ & $ \overline{P}^{(\downarrow)} $ & $ \eta $ \\
	(MWh)     &       (MWh)      &              (MW)            &             (MW)               &          \\ 
	\midrule
	0        &         80       &               30              &               30                & 0.81     \\ 
	\bottomrule
	\end{tabular}
\end{table}
									
\subsection{VPP with wind farm}\label{subsec:VPP_wind}
					
\begin{figure}[t!]
	\centering
	\includegraphics[width=0.8\columnwidth]{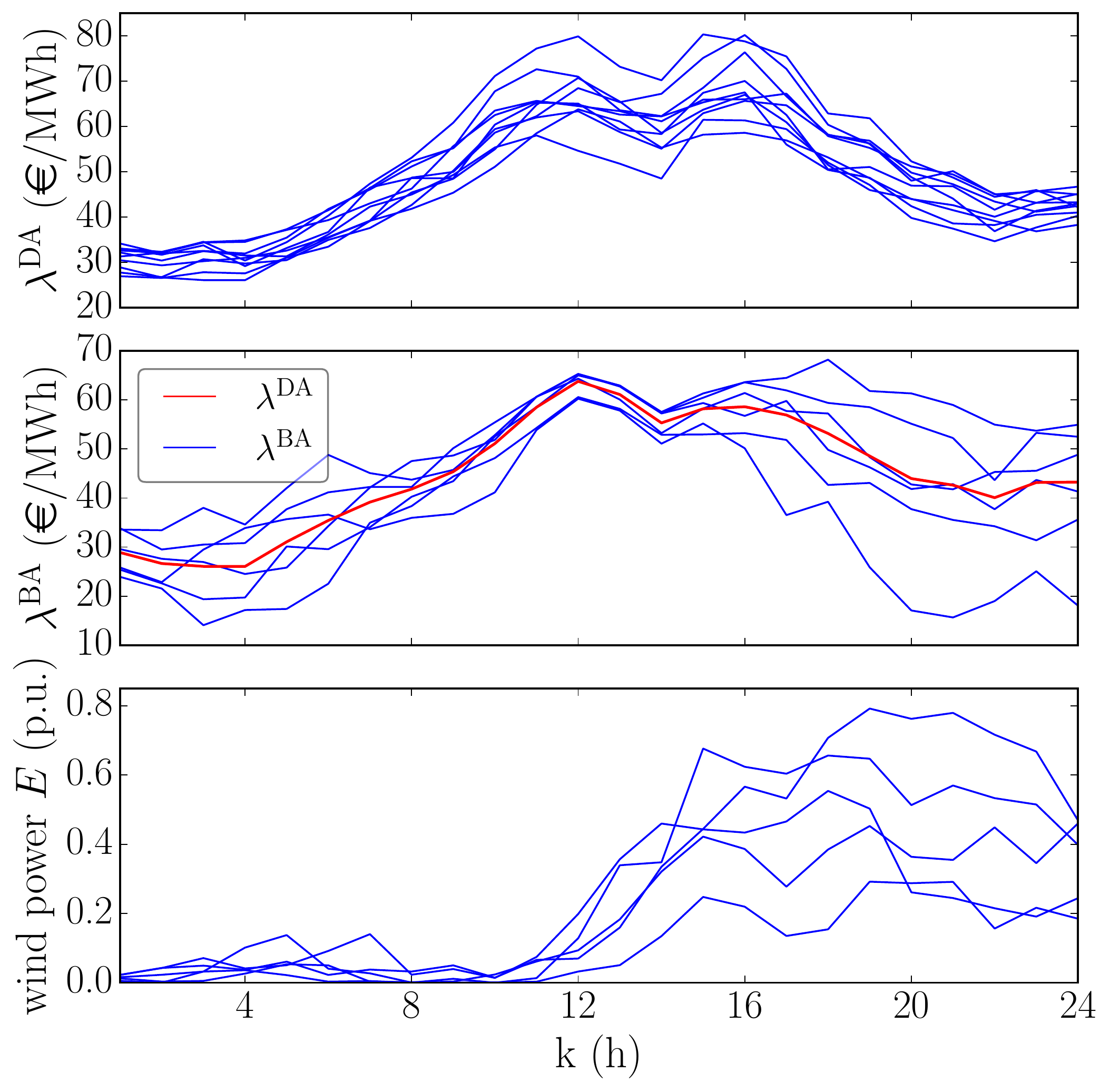}
	\caption{Input day-ahead market price (top), balancing price (middle), and wind power production scenarios (bottom) to the offering model.} \label{fig:Scenarios_wind}
\end{figure}

	Fig.~\ref{fig:Scenarios_wind} shows the ten selected trajectories for the day-ahead market price $\lambda^{\textrm{DA}}_{ik}$, the six chosen balancing price scenarios $\lambda^{\textrm{BA}}_{ijk}$ for a sample day-ahead trajectory $i$, and the five selected trajectories for the wind power production $E_{\omega k}$ (in p.u.).
					
\begin{figure}[b!]
	\centering
	\includegraphics[width=0.75\columnwidth]{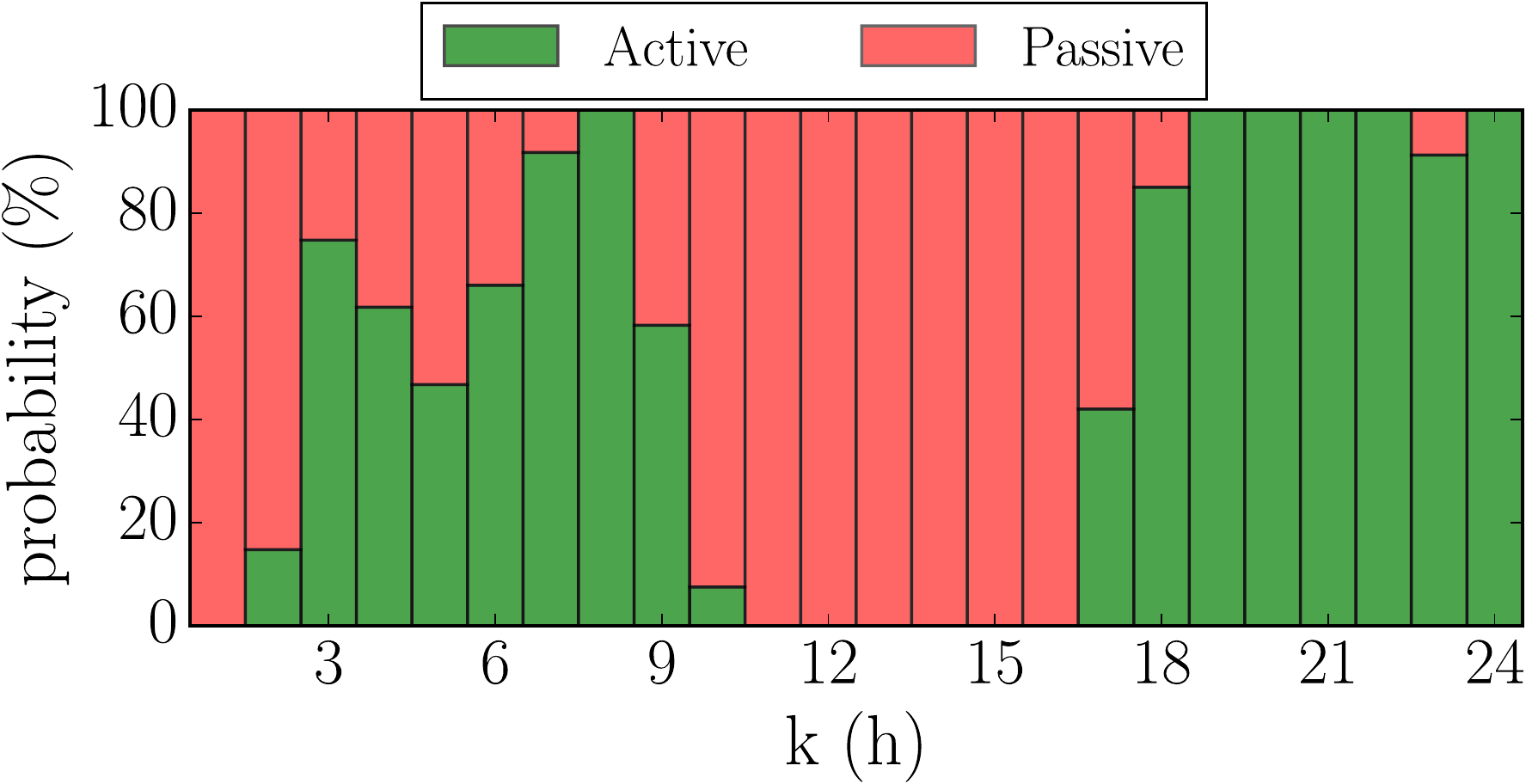}
	\caption{Probability of being active vs. passive ($\overline{E}=50$ MW).} \label{fig:Act_Pas_probs_Wind}
\end{figure}

	The wind farm capacity $\overline{E}$ is initially set to 50 MW. 
    We solve the \textit{Active/Passive} offering model \eqref{eq:Off_Strategy} using as input the scenarios shown in Fig.~\ref{fig:Scenarios_wind}. 
    The complementarity between the active/passive choice is enforced through the binary variables $\epsilon_{ik}$. If $\epsilon_{ik}=1$, the VPP is predicting to act as an active participant during the interval $k$ of the balancing stage, provided that the day-ahead price scenario $i$ realizes. For the same scenario $i$ and interval $k$, if $\epsilon_{ik}=0$, then the VPP is expecting to behave passively. Being $\epsilon^*_{ik}$ the optimal solution, then the probability that the VPP will be active is computed as $\sum_i \pi^{\mathrm{DA}}_i \epsilon^*_{ik}$, and the probability that it will be passive as $\sum_i \pi^{\mathrm{DA}}_i(1- \epsilon^*_{ik})$. These probabilities are illustrated in Fig.~\ref{fig:Act_Pas_probs_Wind}. 
					
\begin{figure}[t!]
	\centering
	\includegraphics[width=0.85\columnwidth]{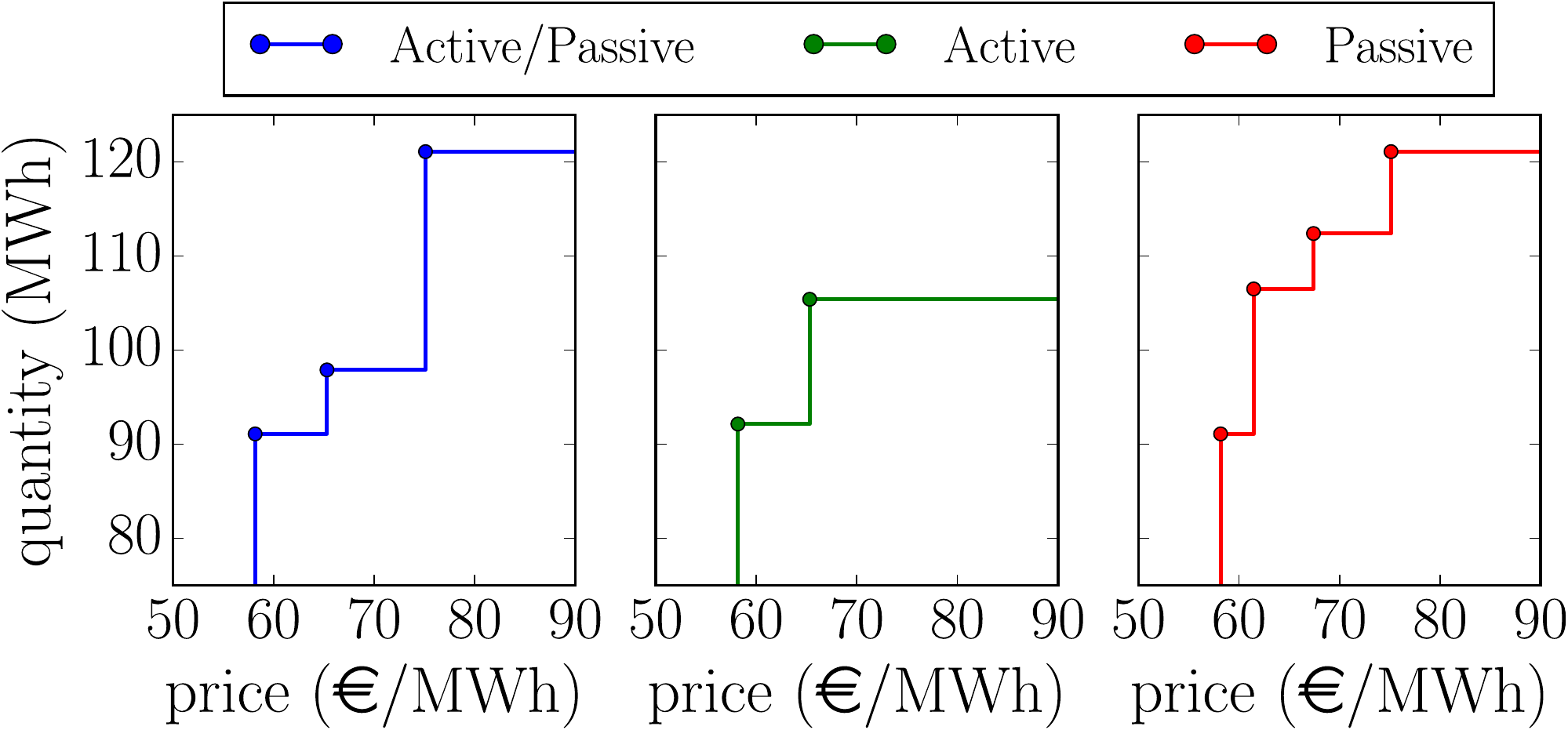}
	\caption{Day-ahead market offer curves at $k=15$ ($\overline{E}=50$ MW).} \label{fig:DA_curves_wind}
\end{figure}

	From midnight to 10 a.m., the VPP will decide to be either active or passive depending on the day-ahead price realization $i$. On one hand, the uncertainty in the wind power production is limited in this trading interval, which would benefit an active approach as the flexibility of the controllable units could be used to offer regulating energy. On the other hand, even if uncertain, the spread between the balancing price scenarios and the day-ahead prices is also limited, leading consequently to low additional profits resulting from an active participation. Then, from 10 a.m. to 4 p.m., the VPP decides to be a passive participant for each realization $i$ of the day-ahead price. Indeed, the uncertainty in the balancing market prices is limited while the amount of wind power production is very uncertain. Finally, from 6 p.m. to midnight, the VPP is almost sure to sell regulating energy in the balancing market. This translates into internally handling the forecasting error of wind power production, which is highly uncertain in these time intervals as it can vary from 20\% to almost 80\% of the wind farm capacity. However, the balancing market price scenarios are going to be far-off the day-ahead market price with high probability; accordingly, passive deviations from the day-ahead schedule may result in heavy penalties, while selling regulating energy can be very profitable. 
\\
\indent
As an example, Fig.~\ref{fig:DA_curves_wind} illustrates the day-ahead market offer curves for the interval $k=15$ obtained when using the \textit{Active/Passive} strategy (left), \textit{Active} strategy (middle), and \textit{Passive} strategy (right). Note that the output of offering model \eqref{eq:Off_Strategy} is a discrete set of price-quantity pairs ($q^{\mathrm{DA}}_{ik},\lambda^{\mathrm{DA}}_{ik}$). To obtain the step functions of Fig.~\ref{fig:DA_curves_wind} starting from the set of price-quantity pairs, we use the standard methodology presented in \cite{conejo2010decision} and \cite{mazzi2017price}. 
    From Fig.~\ref{fig:DA_curves_wind}, we note that the  \textit{Active} approach appears to be less ``reactive" to the day-ahead market price compared to the other strategies, and no additional quantity is scheduled for high values of the day-ahead market price. Indeed, the position of the VPP after the day-ahead market affects its capability to internally compensate for the wind power uncertainty. Therefore, the VPP position is more constrained and driven by feasibility limitations compared to the other two strategies.					
\\
\indent              
Finally, in Table \ref{tab:VPP_Erho_wind} we compare the expected profit of the three offering models for different values of the wind farm capacity $\overline{E}$, ranging from 10 to 90 MW. When $\overline{E}=10$ MW, the expected profit increment under the \textit{Active/Passive} approach is 0.6\% and 8.1\% compared to the \textit{Active} and \textit{Passive} strategies, respectively. In effect, if the capacity of the stochastic unit is small, then the VPP can internally handle most of the wind power deviations and offer its regulating energy into the balancing market. Accordingly, the increase in profit compared to the \textit{Active} strategy is limited whereas  the \textit{Passive} strategy is strongly outperformed. This trend progressively changes as the wind farm capacity $\overline{E}$ increases. As $\overline{E}$ grows, the VPP based on an \textit{Active/Passive} participation is more likely to settle deviations in the balancing stage and has less flexibility to offer in the balancing market. When $\overline{E}=90$ MW, the increase in profit is 4.6\% compared to the \textit{Active} strategy and 2.1\% compared to the \textit{Passive} one. 
                    
\begin{table}[h!]
	\centering
	\caption{Expected profit for the three offering strategies for different values of the wind farm capacity $\overline{E}$.}
	\label{tab:VPP_Erho_wind}
	\begin{tabular}{cccc}
		\toprule
		$\overline{E}$ (MW)   & \multicolumn{3}{c}{Expected profit ($10^3$ \euro)  }  \\
		& \textit{Active/Passive} & \textit{Active} & \textit{Passive} \\ 
		\midrule                                                                                                                               10             &  18.16                  & 18.06           & 16.81            \\
		30             &  22.89                  & 22.48           & 21.69            \\
		50             &  27.59                  & 26.85           & 26.57            \\
		70             &  32.35                  & 31.16           & 31.45            \\
		90             &  37.07                  & 35.44           & 36.32            \\ 
		\bottomrule
	\end{tabular}
\end{table}
					
\subsection{VPP with PV solar}\label{subsec:VPP_solar}
					
\begin{figure}[b!]
	\centering
	\includegraphics[width=0.8\columnwidth]{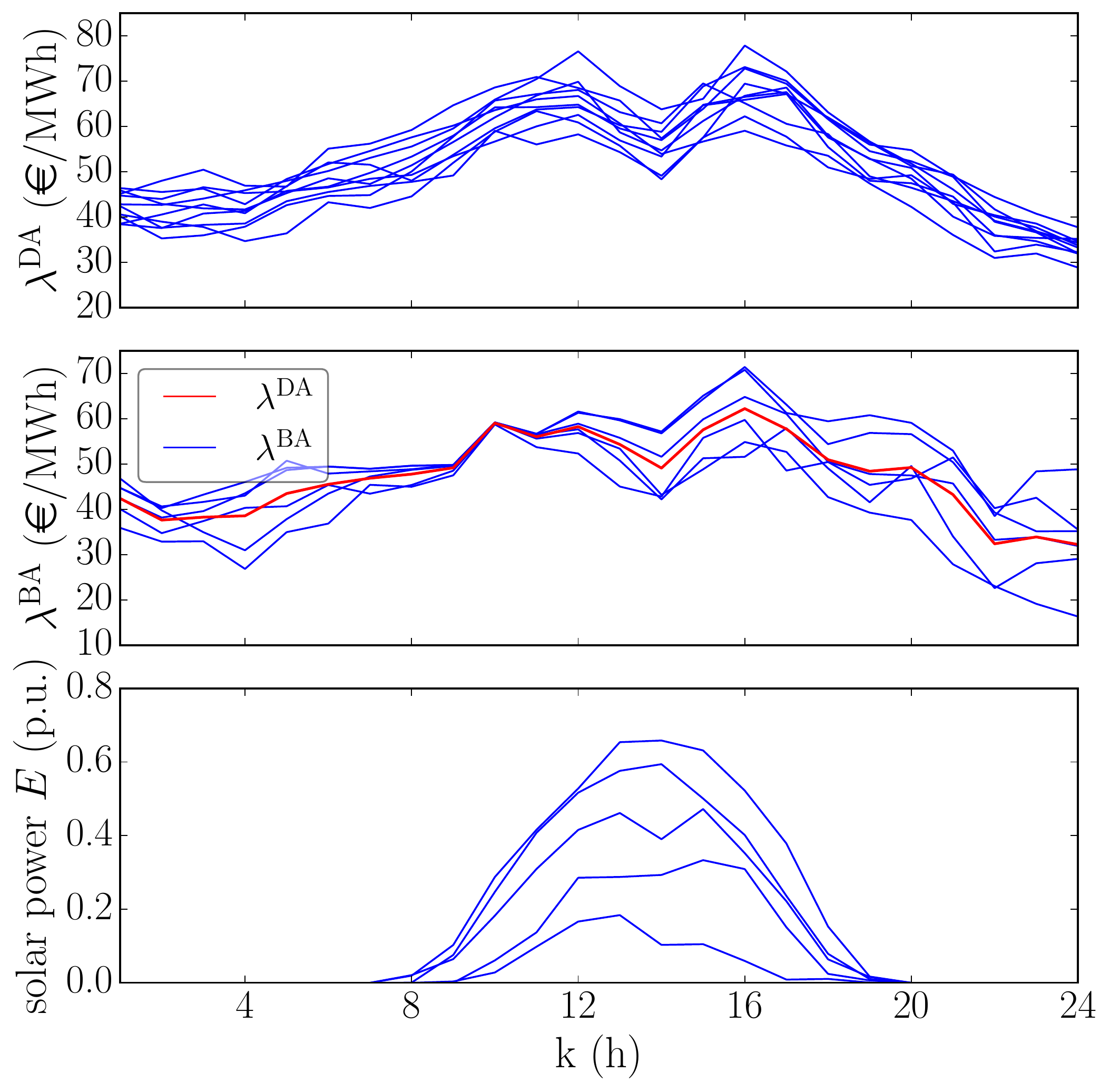}
	\caption{Input day-ahead market price (top), balancing price (middle), and PV solar power production scenarios (bottom) to the offering model.}\label{fig:Scenarios_solar}
\end{figure}

	Fig.~\ref{fig:Scenarios_solar} shows the ten selected trajectories for the day-ahead market price $\lambda^{\textrm{DA}}_{ik}$, the six chosen balancing price scenarios $\lambda^{\textrm{BA}}_{ijk}$ for a sample day-ahead trajectory $i$, and the five selected trajectories for the solar power production $E_{\omega k}$ (in p.u.).
					
\begin{figure}[b!]
		\centering
		\includegraphics[width=0.75\columnwidth]{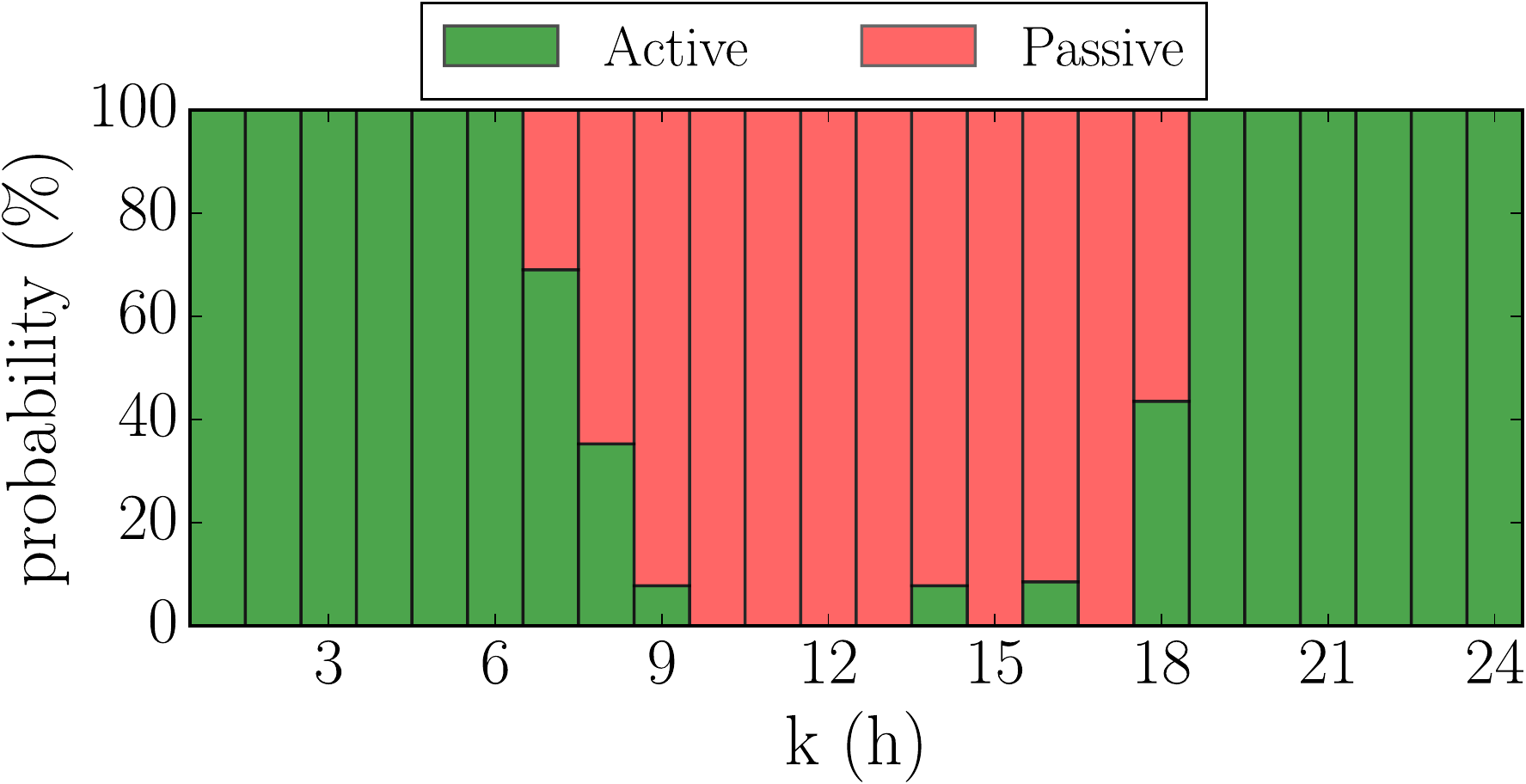}
		\caption{Probability of being active vs. passive ($\overline{E}=50$ MW).} \label{fig:Act_Pas_probs_Solar}
	\end{figure}

	The PV solar unit is initially considered with capacity $\overline{E}=50$ MW. We run the \textit{Active/Passive} offering model \eqref{eq:Off_Strategy}, and in Fig.~\ref{fig:Act_Pas_probs_Solar} show the probabilities of the VPP being active and passive in the trading horizon, computed as in Section \ref{subsec:VPP_wind}.
From midnight to 6 a.m. and from 8 p.m. to midnight, the VPP decides to be active in the balancing market for each day-ahead scenario $i$. Indeed, these time intervals are before and after the sunset, and the VPP is certain that the output of its PV solar unit will be zero. Differently, from 10 a.m. to 5 p.m., the VPP is almost sure that it will passively deviate from its contracted schedule to compensate for the forecasting error of the PV solar unit. In this time horizon, the PV power production is very uncertain (e.g., at 2 p.m. it can vary from 20\% to 70\% of the unit capacity), and it is more convenient to settle deviations in the balancing market. Finally, from 6 a.m. to 9 a.m., the VPP will decide whether the active or the passive participation is more profitable depending on the day-ahead price realization $i$, and the associated amount of energy contracted. In this time interval, the uncertainty of the PV solar production is in fact limited, which would suggest that an active participation may be preferable. However, the possibility of gaining extra profits from the balancing market is low as the balancing market price scenarios are very close to the day-ahead price. Differently, from 6 p.m. to 8 p.m., an active participation is more attractive since the balancing market price scenarios give the opportunity to gain extra profits.
					
\begin{figure}[t!]
	\centering
	\includegraphics[width=0.85\columnwidth]{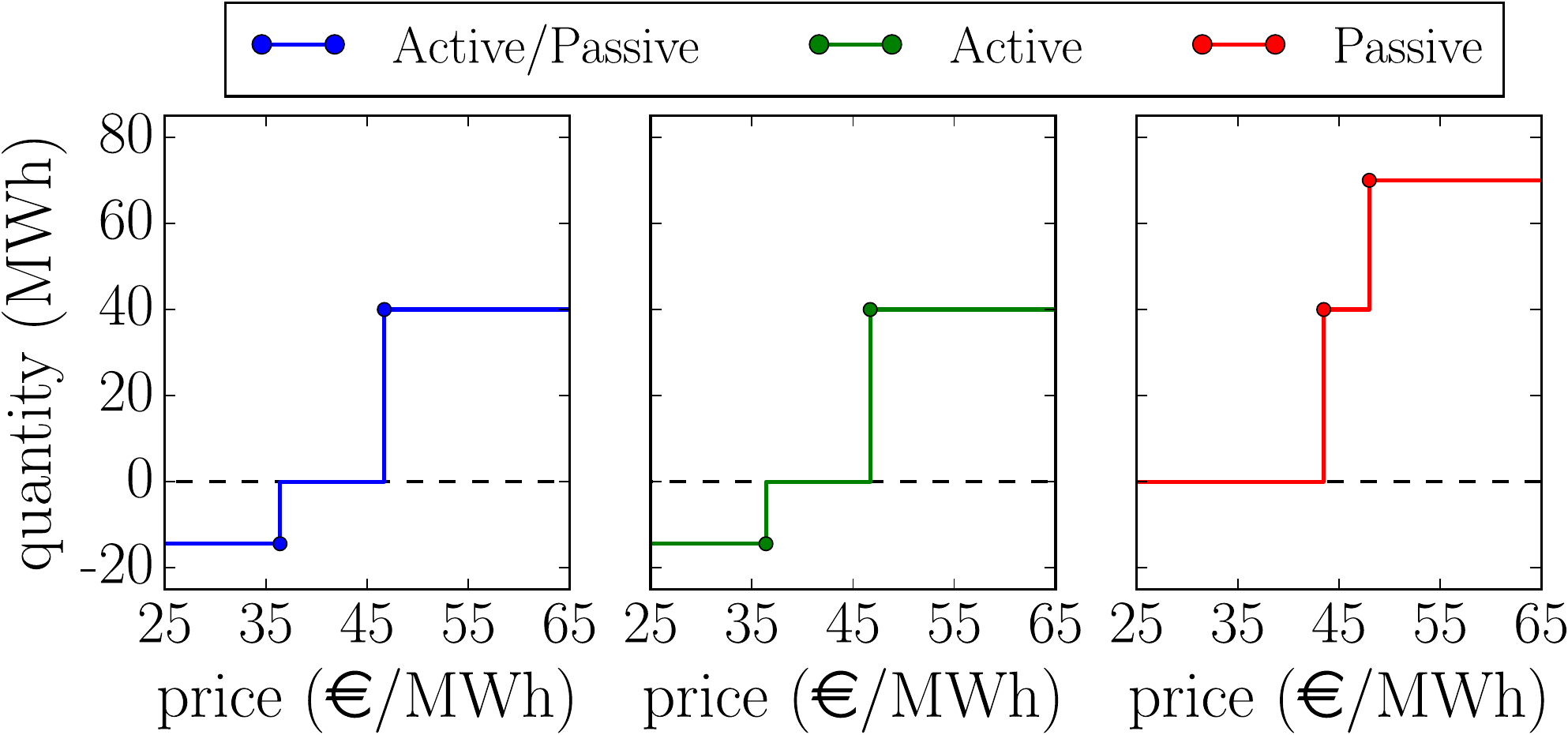}
	\caption{Day-ahead market offer curves at $k=5$. ($\overline{E}=50$ MW).} 	\label{fig:DA_curves_solar}
\end{figure}
	In Fig.~\ref{fig:DA_curves_solar}, we plot the day-ahead market offer curves from the three VPP strategies \textit{Active/Passive} (left), \textit{Active} (middle), and \textit{Passive} (right) at $k=5$. In this interval, the offer curves derived for the \textit{Active/Passive} and \textit{Active} strategy are equivalent, which is consistent with the results shown in Fig.~\ref{fig:Act_Pas_probs_Solar} where the VPP decides to be always active from midnight to 6 a.m.. 
	Compared to these two models, the \textit{Passive} one scheduled more energy in the day-ahead as the VPP is less driven by feasibility constraints. 
					
	Lastly, we compare the expected profit obtained with the three offering models for increasing capacity of the PV unit, $\overline{E}$, from 10 to 90 MW. We show the results in Table \ref{tab:VPP_Erho_solar}. 
\begin{table}[b!]
	\centering
	\caption{Expected profit for the three offering strategies for different values of the PV solar unit capacity $\overline{E}$.}
	\label{tab:VPP_Erho_solar}
	\begin{tabular}{cccc}
		\toprule
		$\overline{E}$ (MW)   & \multicolumn{3}{c}{Expected profit ($10^3$ \euro)  }  \\
		& \textit{Active/Passive}   & \textit{Active}  & \textit{Passive} \\  
		\midrule                                                                                                                               10             &  18.59                    & 18.27            & 17.54            \\
		30             &  22.19                    & 21.28            & 21.15            \\
		50             &  25.79                    & 24.24            & 24.76            \\
		70             &  29.39                    & 27.17            & 28.36            \\
		90             &  32.98                    & 30.08            & 31.97            \\ 
		\bottomrule
	\end{tabular}
\end{table}
When $\overline{E}=10$ MW, the expected profit of the \textit{Active/Passive} and the \textit{Active} approach are similar. Indeed, as the capacity of the PV unit is small, the VPP can easily handle the PV solar uncertainty internally and offer the remaining flexibility in the balancing market. The expected profit increase under the \textit{Active/Passive} approach is 1.7\% and 6.0\% compared to the \textit{Active} and the \textit{Passive} strategies, respectively. As the PV unit capacity increases, the \textit{Passive} strategy becomes more competitive as it can benefit from a more flexible VPP operation, and can consequently contract more profitable positions in the day-ahead market. Instead, as $\overline{E}$ grows, the \textit{Active} approach becomes increasingly constrained in its operation. First, it has less flexibility to offer in the balancing stage as it needs to allocate it to balance the PV unit forecasting errors. Second, the day-ahead position is more constrained by ensuring a feasible real-time operation, thus it is less driven by the market prices.
					
	\section{Conclusions}\label{sec:Conclusions}
					
	In this paper we proposed an innovative participation model for the balancing market, denoted by \textit{Active/Passive}, aimed to increase the flexibility of market participants as well as the amount of regulating energy available to the TSO in the real-time. Specifically, we suggest to allow market agents such as VPPs to actively offer regulating energy in time intervals where they can ensure to internally handle the eventual forecasting errors of the stochastic energy sources, while passively deviating from their day-ahead schedule in other intervals. We enforced these two participation modes (active and passive) to be complementary, and agents submitting regulating energy offers for a specific trading interval are prevented from creating an imbalance in the same interval.
					
	To analyze this novel participation model, we took the perspective of a VPP that includes both controllable and stochastic generation units, and that trades in a two-settlement electricity market. The \textit{Active/Passive} offering strategy arises as a three-stage decision making problem. We formulated this problem as a MILP in which binary variables are introduced to model the feasible operating region of conventional production units and the complementarity between active and passive participation in the balancing stage. Compared to an \textit{Active} and a \textit{Passive} strategy, computational experiments showed that an \textit{Active/Passive} approach can result in a significantly higher VPP expected profit (up to 8\% higher). The analysis reveals that the active participation is more attractive for the VPP in the hourly intervals with limited production uncertainty from the stochastic sources and profitable balancing market price scenarios, and the passive one when highly uncertain renewable energy production is combined with narrow balancing market price scenarios (i.e., close to the day-ahead price). 
					
	The proposed framework is potentially relevant from the perspective of a system operator, who would benefit from having more regulating energy available in real time, and can be also seen as a lever to facilitate the integration of renewable power sources through their aggregation into VPPs.
					
	The focus of this work was to provide useful insights on the proposed \textit{Active/Passive} participation model. Therefore, we presented a case study with 300 scenarios and a simplified operating region of the dispatchable generators to keep the model intuitive and solvable in about 30 minutes. Further research may be in the direction of developing more efficient algorithms capable of solving the model for a larger number of scenarios or including more operating details of the units.

	\ifCLASSOPTIONcaptionsoff
	\newpage
	\fi
	
	\bibliographystyle{IEEEtran}
	\bibliography{bibliography}

\begin{thebibliography}{10}
\providecommand{\url}[1]{#1}
\csname url@samestyle\endcsname
\providecommand{\newblock}{\relax}
\providecommand{\bibinfo}[2]{#2}
\providecommand{\BIBentrySTDinterwordspacing}{\spaceskip=0pt\relax}
\providecommand{\BIBentryALTinterwordstretchfactor}{4}
\providecommand{\BIBentryALTinterwordspacing}{\spaceskip=\fontdimen2\font plus
\BIBentryALTinterwordstretchfactor\fontdimen3\font minus
  \fontdimen4\font\relax}
\providecommand{\BIBforeignlanguage}[2]{{%
\expandafter\ifx\csname l@#1\endcsname\relax
\typeout{** WARNING: IEEEtran.bst: No hyphenation pattern has been}%
\typeout{** loaded for the language `#1'. Using the pattern for}%
\typeout{** the default language instead.}%
\else
\language=\csname l@#1\endcsname
\fi
#2}}
\providecommand{\BIBdecl}{\relax}
\BIBdecl

\bibitem{morales2013integrating}
J.~M. Morales, A.~J. Conejo, H.~Madsen, P.~Pinson, and M.~Zugno,
  \emph{Integrating renewables in electricity markets: operational
  problems}.\hskip 1em plus 0.5em minus 0.4em\relax Springer Science \&
  Business Media, 2013, vol. 205.

\bibitem{bremnes2004}
J.~B. Bremnes, ``Probabilistic wind power forecasts using local quantile
  regression,'' \emph{Wind Energy}, vol.~7, no.~1, pp. 47--54, 2004.

\bibitem{pinson2007trading}
P.~Pinson, C.~Chevallier, and G.~N. Kariniotakis, ``Trading wind generation
  from short-term probabilistic forecasts of wind power,'' \emph{IEEE
  Transactions on Power Systems}, vol.~22, no.~3, pp. 1148--1156, 2007.

\bibitem{bitar2012bringing}
E.~Y. Bitar, R.~Rajagopal, P.~P. Khargonekar, K.~Poolla, and P.~Varaiya,
  ``Bringing wind energy to market,'' \emph{IEEE Transactions on Power
  Systems}, vol.~27, no.~3, pp. 1225--1235, 2012.

\bibitem{dent2011opportunity}
C.~J. Dent, J.~W. Bialek, and B.~F. Hobbs, ``Opportunity cost bidding by wind
  generators in forward markets: Analytical results,'' \emph{IEEE Transactions
  on Power Systems}, vol.~26, no.~3, pp. 1600--1608, 2011.

\bibitem{morales2010short}
J.~M. Morales, A.~J. Conejo, and J.~P{\'e}rez-Ruiz, ``Short-term trading for a
  wind power producer,'' \emph{IEEE Transactions on Power Systems}, vol.~25,
  no.~1, pp. 554--564, 2010.

\bibitem{Arroyo2000}
J.~M. Arroyo and A.~J. Conejo, ``Optimal response of a thermal unit to an
  electricity spot market,'' \emph{IEEE Trans. Power Syst.}, vol.~15, no.~3,
  pp. 1098--1104, 2000.

\bibitem{Arroyo2004}
------, ``Modeling of start-up and shut-down power trajectories of thermal
  units,'' \emph{IEEE Trans. Power Syst.}, vol.~19, no.~3, pp. 1562--1568,
  2004.

\bibitem{Conejo2004}
A.~J. Conejo, F.~J. Nogales, J.~M. Arroyo, and R.~Garc{\'\i}a-Bertrand,
  ``Risk-constrained self-scheduling of a thermal power producer,'' \emph{IEEE
  Trans. Power Syst.}, vol.~19, no.~3, pp. 1569--1574, 2004.

\bibitem{Conejo2002}
A.~J. Conejo, F.~J. Nogales, and J.~M. Arroyo, ``Price-taker bidding strategy
  under price uncertainty,'' \emph{IEEE Trans. Power Syst.}, vol.~17, no.~4,
  pp. 1081--1088, 2002.

\bibitem{Baillo2004}
A.~Baillo, M.~Ventosa, M.~Rivier, and A.~Ramos, ``Optimal offering strategies
  for generation companies operating in electricity spot markets,'' \emph{IEEE
  Trans. Power Syst.}, vol.~19, no.~2, pp. 745--753, 2004.

\bibitem{conejo2010decision}
A.~J. Conejo, M.~Carri{\'o}n, and J.~M. Morales, \emph{Decision Making Under
  Uncertainty in Electricity Markets}.\hskip 1em plus 0.5em minus 0.4em\relax
  Springer, 2010.

\bibitem{Maenhoudt2014}
M.~Maenhoudt and G.~Deconinck, ``Strategic offering to maximize day-ahead
  profit by hedging against an infeasible market clearing result,'' \emph{IEEE
  Trans. Power Syst.}, vol.~29, no.~2, pp. 854--862, 2014.

\bibitem{ruiz2009direct}
N.~Ruiz, I.~Cobelo, and J.~Oyarzabal, ``A direct load control model for virtual
  power plant management,'' \emph{IEEE Transactions on Power Systems}, vol.~24,
  no.~2, pp. 959--966, 2009.

\bibitem{hellmers2016operational}
A.~Hellmers, M.~Zugno, A.~Skajaa, and J.~M. Morales, ``Operational strategies
  for a portfolio of wind farms and chp plants in a two-price balancing
  market,'' \emph{IEEE Trans. Power Syst}, vol.~31, no.~3, pp. 2182--2191,
  2016.

\bibitem{zapata2014comparative}
J.~Zapata, J.~Vandewalle, and W.~D'haeseleer, ``A comparative study of
  imbalance reduction strategies for virtual power plant operation,''
  \emph{Applied Thermal Engineering}, vol.~71, no.~2, pp. 847--857, 2014.

\bibitem{mashhour2011biddingI}
E.~Mashhour and S.~M. Moghaddas-Tafreshi, ``Bidding strategy of virtual power
  plant for participating in energy and spinning reserve markets—part i:
  Problem formulation,'' \emph{IEEE Transactions on Power Systems}, vol.~26,
  no.~2, pp. 949--956, 2011.

\bibitem{mashhour2011biddingII}
------, ``Bidding strategy of virtual power plant for participating in energy
  and spinning reserve markets—part ii: Numerical analysis,'' \emph{IEEE
  Transactions on Power Systems}, vol.~26, no.~2, pp. 957--964, 2011.

\bibitem{peik2013decision}
M.~Peik-Herfeh, H.~Seifi, and M.~Sheikh-El-Eslami, ``Decision making of a
  virtual power plant under uncertainties for bidding in a day-ahead market
  using point estimate method,'' \emph{International Journal of Electrical
  Power \& Energy Systems}, vol.~44, no.~1, pp. 88--98, 2013.

\bibitem{pandvzic2013virtual}
H.~Pand{\v{z}}i{\'c}, I.~Kuzle, and T.~Capuder, ``Virtual power plant mid-term
  dispatch optimization,'' \emph{Applied Energy}, vol. 101, pp. 134--141, 2013.

\bibitem{pandvzic2013offering}
H.~Pand{\v{z}}i{\'c}, J.~M. Morales, A.~J. Conejo, and I.~Kuzle, ``Offering
  model for a virtual power plant based on stochastic programming,''
  \emph{Applied Energy}, vol. 105, pp. 282--292, 2013.

\bibitem{kardakos2016optimal}
E.~G. Kardakos, C.~K. Simoglou, and A.~G. Bakirtzis, ``Optimal offering
  strategy of a virtual power plant: A stochastic bi-level approach,''
  \emph{IEEE Transactions on Smart Grid}, vol.~7, no.~2, pp. 794--806, 2016.

\bibitem{skajaa2015intraday}
A.~Skajaa, K.~Edlund, and J.~M. Morales, ``Intraday trading of wind energy,''
  \emph{IEEE Transactions on Power Systems}, vol.~30, no.~6, pp. 3181--3189,
  2015.

\bibitem{wang2015review}
Q.~Wang, C.~Zhang, Y.~Ding, G.~Xydis, J.~Wang, and J.~{\O}stergaard, ``Review
  of real-time electricity markets for integrating distributed energy resources
  and demand response,'' \emph{Applied Energy}, vol. 138, pp. 695--706, 2015.

\bibitem{mazzi2017price}
N.~Mazzi, J.~Kazempour, and P.~Pinson, ``Price-taker offering strategy in
  electricity pay-as-bid markets,'' \emph{IEEE Transactions on Power Systems},
  2017.

\bibitem{pierro2016multi}
M.~Pierro, F.~Bucci, M.~De~Felice, E.~Maggioni, D.~Moser, A.~Perotto, F.~Spada,
  and C.~Cornaro, ``Multi-model ensemble for day ahead prediction of
  photovoltaic power generation,'' \emph{Solar Energy}, vol. 134, pp. 132--146,
  2016.

\bibitem{pinson2009probabilistic}
P.~Pinson, H.~Madsen, H.~A. Nielsen, G.~Papaefthymiou, and B.~Kl{\"o}ckl,
  ``From probabilistic forecasts to statistical scenarios of short-term wind
  power production,'' \emph{Wind Energy}, vol.~12, no.~1, pp. 51--62, 2009.

\bibitem{pinson2012evaluating}
P.~Pinson and R.~Girard, ``Evaluating the quality of scenarios of short-term
  wind power generation,'' \emph{Applied Energy}, vol.~96, pp. 12--20, 2012.

\bibitem{growe2003scenario}
N.~Growe-Kuska, H.~Heitsch, and W.~Romisch, ``Scenario reduction and scenario
  tree construction for power management problems,'' in \emph{IEEE Power Tech
  Conference Proceedings}, 2003, pp. 1--7.

\end{thebibliography}
		
	\begin{IEEEbiographynophoto}{Nicol\`o Mazzi}
		received the B.Sc. degree in 2011 and the M.Sc. degree in 2014, in energy engineering, and the Ph.D. degree in 2018, in industrial engineering, from the University of Padova, Padova, Italy. He is currently a postdoctoral researcher at the School of Mathematics, University of Edinburgh. \newline
		\indent His research interests include power systems planning and operations, electricity markets, stochastic and hierarchical optimization, and decomposition techniques.
	\end{IEEEbiographynophoto}
	
	\begin{IEEEbiographynophoto}{Alessio Trivella}
		received the B.Sc. and M.Sc. degrees in mathematics from the University of Milan in Italy in 2010 and 2012, respectively, and the Ph.D. degree in operations research from the Technical University of Denmark in 2018. He is currently a postdoctoral researcher with the Department of Management Engineering, Technical University of Denmark. \newline
		\indent His research interests include energy operations and investments, sustainable energy integration, stochastic optimization, and approximate dynamic programming.
	\end{IEEEbiographynophoto}
	
	\begin{IEEEbiographynophoto}{Juan M.\ Morales}
        (S'07-M'11-SM'16) received the Ingeniero Industrial degree from the University of  M\'alaga, M\'alaga, Spain, in 2006, and a Ph.D. degree in Electrical Engineering from the University of Castilla-La Mancha, Ciudad Real, Spain, in 2010. He is currently an associate professor in the Department of Applied Mathematics at the University of M\'alaga in Spain. \newline
        \indent His research interests are in the fields of power systems economics, operations and planning; energy analytics and optimization; smart grids; decision-making under uncertainty, and electricity markets.
	\end{IEEEbiographynophoto}
	
	\vfill
	
\end{document}